%% file: YuXuTongJournal.tex
\documentclass[journal]{IEEEtran}
\usepackage{algorithmic}
\usepackage{algorithm}
\usepackage{mathrsfs,amsthm}
\IEEEoverridecommandlockouts
\makeatletter
\def\ps@headings{%
\def\@oddhead{\mbox{}\scriptsize\rightmark \hfil \thepage}%
\def\@evenhead{\scriptsize\thepage \hfil \leftmark\mbox{}}%
\def\@oddfoot{}%
\def\@evenfoot{}}
\makeatother 
\usepackage{graphicx,amssymb,stfloats,subfigure,multicol,amsmath}%
\usepackage{epsfig,epsf,psfrag,amssymb,amsfonts,latexsym,graphicx,mathrsfs,subfigure}
\usepackage{dsfont}
\usepackage{chngpage}
\usepackage[dvips]{color}
\usepackage{textcomp}
\usepackage[normalem]{ulem}
\usepackage{cleveref}
\usepackage{caption}
\newtheorem{Def}{Definition}
\newtheorem{corollary}{Corollary}
\newtheorem{thm}{Theorem}
\newtheorem{lem}{Lemma}
\newtheorem{prop}{Proposition}

 \def\old#1{}    

\input ZYmacros

\begin{document}

\title{Deadline Scheduling as Restless Bandits}
\author{\Large Zhe Yu$^\dagger$, Yunjian Xu$^\ddagger$, and Lang Tong$^\dagger$
\thanks{\scriptsize
Z. Yu$^\dagger$ and L. Tong$^\dagger$ are with the School of Electrical and Computer
Engineering, Cornell University, Ithaca, NY 14853, USA. Y. Xu$^\ddagger$ is with the School of Engineering Systems and Design Pillar, Singapore University of Technology and Design, Singapore, 487372. Email:
{\tt\{zy73,lt35\}@cornell.edu, yunjian\_xu@sutd.edu.sg}. This work is supported in part by the National Science Foundation under Grant CNS-
1248079 and 1549989. An earlier version that contains a subset of the results presented in this paper appeared in  \cite{Yu&Xu&Tong:2016Allerton}.}}

\maketitle

\begin{abstract}
The problem of stochastic deadline scheduling  is considered. A constrained Markov decision process model is introduced in which jobs arrive randomly at a service center with stochastic job sizes, rewards, and completion deadlines.  The service provider faces random processing costs, convex non-completion penalties, and a capacity constraint that limits the simultaneous processing of jobs. Formulated as a restless multi-armed bandit problem, the stochastic deadline scheduling problem is shown to be indexable.  A closed-form expression of the Whittle's index is obtained for the case when the processing costs are constant. An upper bound on the gap-to-optimality for the Whittle's index policy is obtained, and it is shown that  the bound converges to zero as the job arrival rate and the number of available processors increase simultaneously to infinity.
\end{abstract}

\begin{IEEEkeywords}
Constrained Markov decision processes; Restless multi-armed bandits; Stochastic deadline scheduling;
 Whittle's index.
\end{IEEEkeywords}

\section{Introduction}

\input intro_v11

\section{Problem Formulation} \label{sec:II}

\input formulation_v8

\section{Whittle's Index Policy}\label{sec:III}
\input whittleIndex_v7

\section{Performance of Whittle's Index Policy for Finite-Armed Restless Bandits}\label{sec:IV}
\input finiteArm

\section{Asymptotic performance of the Whittle's Index Policy}\label{sec:V}
\input asymptotic_v8
\section{Numerical Results}\label{sec:VI}
\input simulation_v8
\section{Conclusion}\label{sec:VII}
\input conclusion_v1

\begin{appendices}

\input proof_v11
\end{appendices}

{
\bibliographystyle{ieeetran}
\bibliography{Bibs/Journal,Bibs/Conf,Bibs/Book,Bibs/Misc}
}


\vspace{-1em}
\begin{IEEEbiography}[{\includegraphics[width=1in,height=1.25in,clip,keepaspectratio]{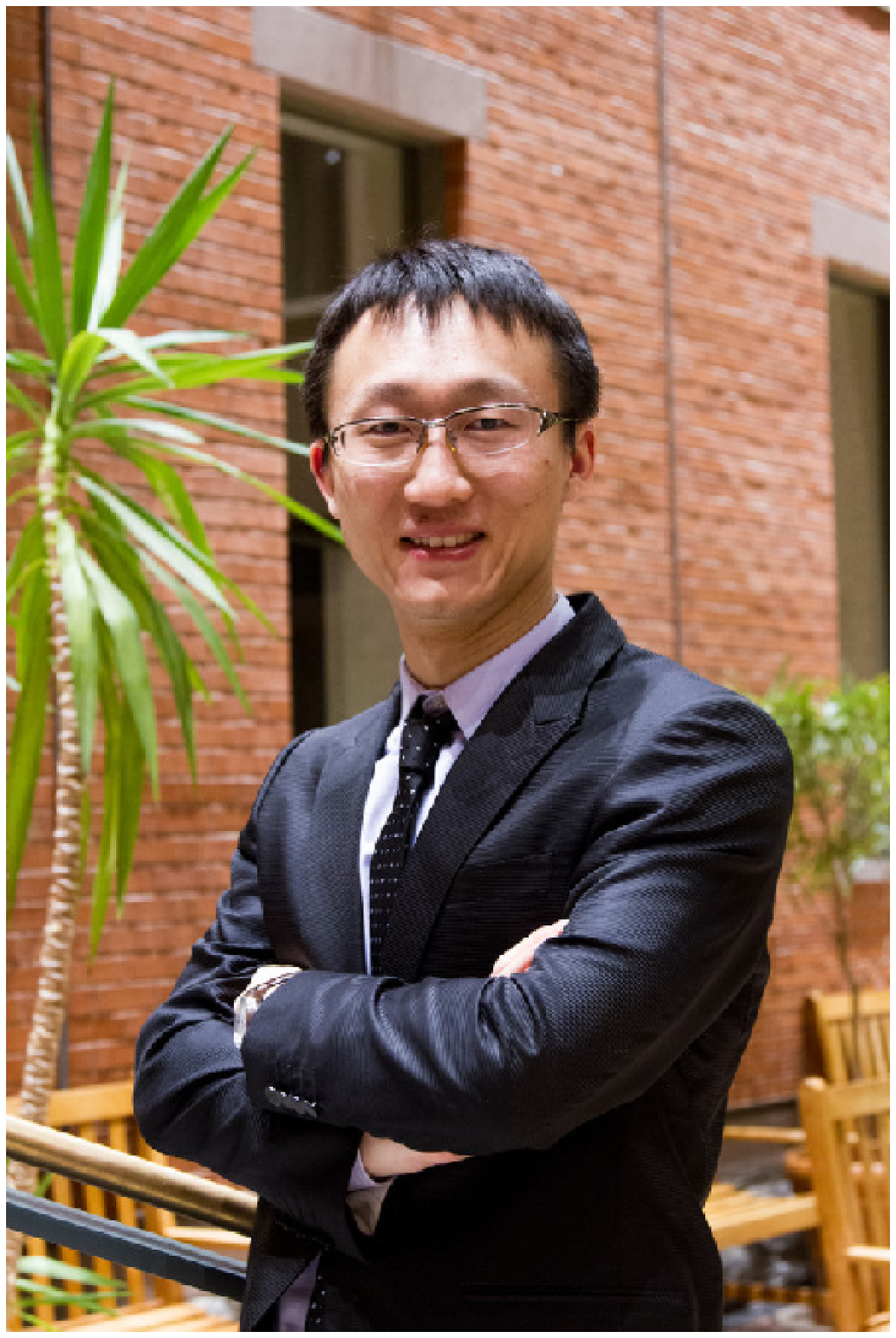}}]{Zhe Yu}
received his B.E. degree from Department
of Electrical Engineering, Tsinghua University, Beijing, China
in 2009, M.S. degree from Department
of Electrical and Computer Engineering, Carnegie Mellon University, Pittsburgh, PA, USA
in 2010, and Ph.D. from the School of Electrical and
Computer Engineering, Cornell University, Ithaca, NY, USA in 2016. He joined State Grid Research Institute North America in 2017. His
current research interests focus on power system and smart grid, demand response, dynamic programming and optimization.
\end{IEEEbiography}
\vspace{-1em}
\begin{IEEEbiography}[{\includegraphics[width=1in,height=1.25in,clip,keepaspectratio]{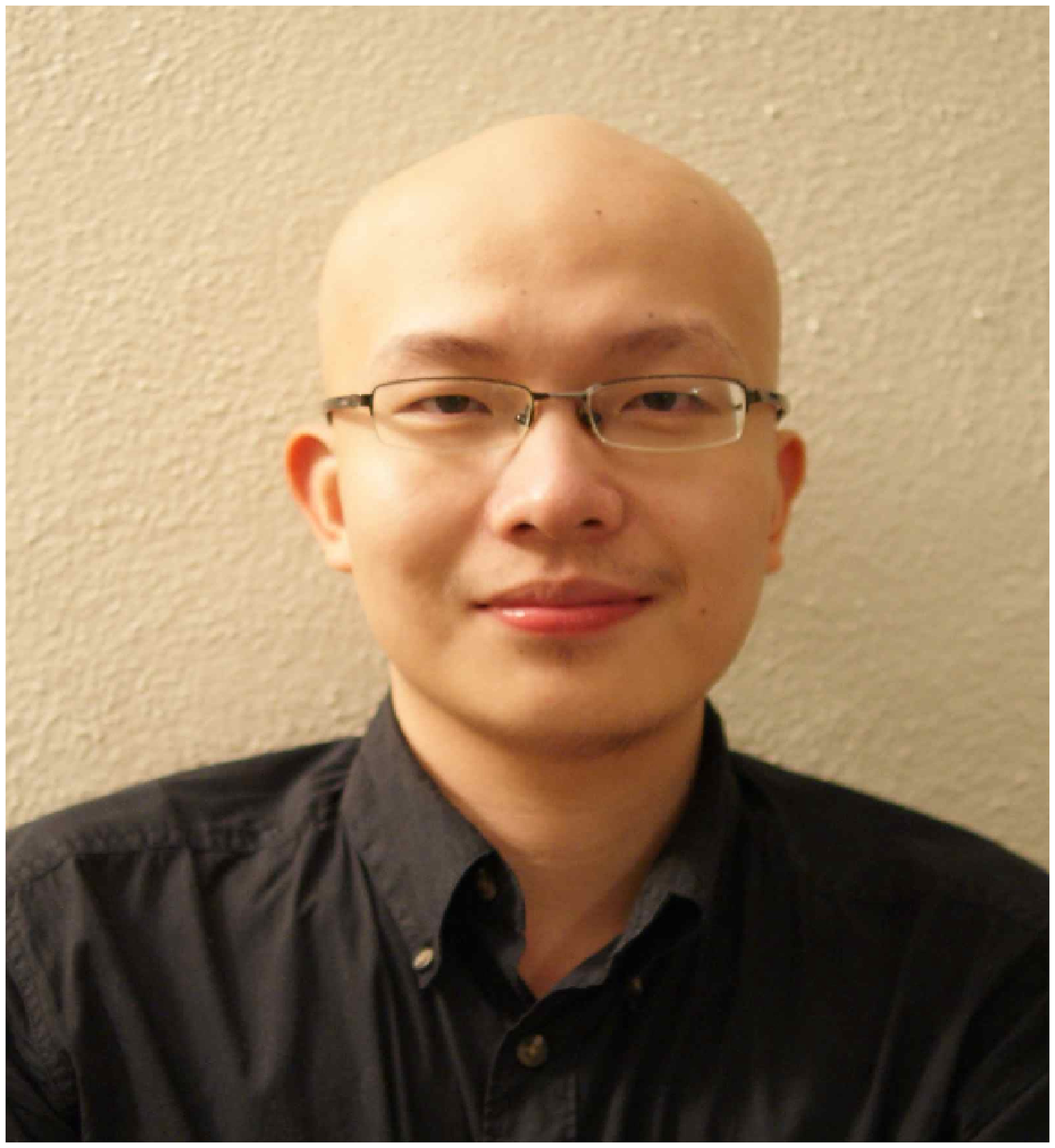}}]{Yunjian Xu}(S'06-M'10)
received the B.S. and M.S. degrees in Electrical
Engineering from Tsinghua University, Beijing, China, in 2006 and
2008, respectively, and the Ph.D. degree from the Massachusetts
Institute of Technology (MIT), Cambridge, MA, USA, in 2012.

Dr. Xu was a CMI (Center for the Mathematics of Information)
postdoctoral fellow at the California Institute of Technology,
Pasadena, CA, USA, in 2012-2013. Before joining the Chinese University
of Hong Kong (CUHK) as an assistant professor, he was an
assistant professor at the Singapore University of
Technology and Design in 2013-2017.
His research interests focus on power systems and electricity markets, with
emphasis on power system control and optimization, wholesale
electricity market design, and the aggregation of distributed energy
resources in power distribution systems.
\end{IEEEbiography}
\vspace{-1em}
\begin{IEEEbiography}[{\includegraphics[width=1in,height=1.25in,clip,keepaspectratio]{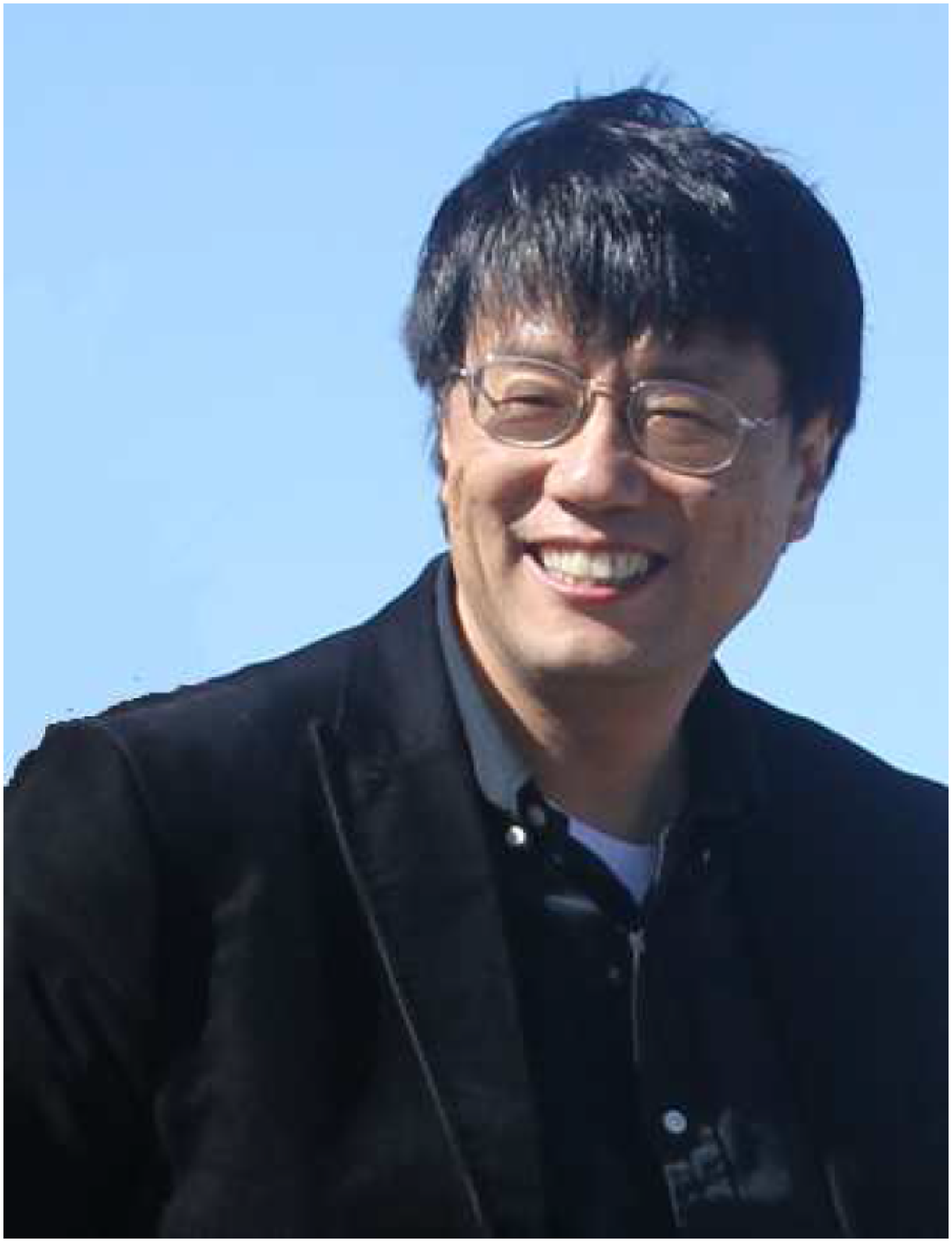}}]{Lang Tong}
(F'05) is the Irwin
and Joan Jacobs Professor in Engineering
of Cornell University and the site director of
Power Systems Engineering Research Center
(PSERC). He received the B.E. degree from
Tsinghua University in 1985, and M.S. and
Ph.D. degrees in electrical engineering in 1987
and 1991, respectively, from the University
of Notre Dame. He was a Postdoctoral Research
Affiliate at the Information Systems
Laboratory, Stanford University in 1991. He
was the 2001 Cor Wit Visiting Professor at the Delft University of
Technology and had held visiting positions at Stanford University
and the University of California at Berkeley.
Lang Tong's research is in the general area of statistical inference,
communications, and complex networks. His current research focuses
on inference, optimization, and economic problems in energy and
power systems. He received the 1993 Outstanding Young Author
Award from the IEEE Circuits and Systems Society, the 2004 best
paper award from IEEE Signal Processing Society, and the 2004
Leonard G. Abraham Prize Paper Award from the IEEE Communications
Society. He is also a coauthor of seven student paper
awards. He received Young Investigator Award from the Office of
Naval Research. He was a Distinguished Lecturer of the IEEE Signal
Processing Society.
\end{IEEEbiography}

\end{document}

%% file: ZYmacros.tex
\def\beq{\begin{equation}}
\def\eeq{\end{equation}}
\def\bea{\begin{eqnarray}}
\def\eea{\end{eqnarray}}
\def\ba{\begin{array}}
\def\ea{\end{array}}
\def\bitem{\begin{itemize}}
\def\eitem{\end{itemize}}
\def\benum{\begin{enumerate}}
\def\eenum{\end{enumerate}}
\def\defeq{{\stackrel{\Delta}{=}}}
\def\edoc{\end{document}}

\def\eg{{\it e.g., \/}}

\def\ie{{\it i.e.,\ \/}}

\definecolor{bgrd}{rgb}{1,1,1}
\definecolor{grey}{rgb}{0.9,0.9,0.6}
\definecolor{gray}{rgb}{0.5,0.5,0.5}
\definecolor{dkr}{rgb}{0.6,0.2,0.2}
\definecolor{dkg}{rgb}{0,0.5,0}
\definecolor{dkb}{rgb}{0.0,0.1,0.7}
\definecolor{light-gray}{gray}{0.85}

\newcommand{\Smsc}{\mathscr{S}}

\newcommand{\Xmsc}{\mathscr{X}}

\def\abf{{\bf a}}

%% file: intro_v11.tex
The deadline scheduling problem, in its most generic setting, is the scheduling of jobs with different workloads and deadlines for completion.  Typically, not enough servers are available to satisfy all the demand; the cost of processing may vary with time, and unfinished jobs incur penalties.

In this paper, we are interested in the {\em stochastic deadline scheduling problem} where key parameters of the problem such as job arrivals, workloads, deadlines of completion, and processing costs are stochastic. In particular, we consider  the problem of maximizing the discounted rewards over an infinite scheduling horizon. 

A prototype application of such a problem is the charging of electric vehicles (EVs) at a charging service center \cite{Yu&Chen&Tong:16CSEE,Yu&Xu&Tong:15Allerton}.  In such applications,  EVs arrive at the service center randomly, each with its own charging demand and deadline for completion.  The charging cost depends on the cost of electricity at the time of charging, and a penalty is imposed when the service provider is unable to fulfill the request.  
Similar applications include the scheduling of packet transmission for real-time wireless networks \cite{Hou&Kumar:13SLCN}, of jobs at data centers \cite{Vilaplana&Etal:14JS}, of nursing personnel in hospitals \cite{Waner&Prawda:MS1972}, for internet streaming \cite{Chen&Primet:07CCG},  and at customer service centers \cite{Dai&He:11TOR}.

The stochastic deadline scheduling problem is an instance of stochastic dynamic programming, for which obtaining the optimal solution is fundamentally intractable.  However, practical applications often mandate that the processing schedule be constructed in real time. This means that, in general, one may have to sacrifice optimality in favor of approximate solutions that are scalable algorithmically and have performance close to that of the optimal scheduler.
An important class of such algorithms is the so-called {\it index policies} \cite{Gittins:79JRSS} that attach an index to each unfinished job, rank them according to their indices, and assign available processors to the top ranked jobs. The index of each job is determined by the state of the job itself and independent of the states of other jobs. Such policies offer scalable solutions if the index and ranking algorithm can be computed online. An index policy becomes especially attractive if its gap-to-optimality can be bounded and shown to be diminishing in cases of practical interest.

\subsection{Summary of Results}
In this paper, we formulate the stochastic deadline scheduling problem as a restless multi-armed bandit (RMAB) problem originally introduced by Whittle \cite{Whittle:1988JAP}.  We  examine the indexability of the problem and the performance of the Whittle's index policy.  To this end, we introduce a constrained Markov decision process (MDP) model with the objective of maximizing the expected (discounted) profit subject to a constraint on the maximum number of jobs that can be processed simultaneously.
The constructed MDP model captures the randomness in job arrivals, job sizes,  deadlines, and processing costs.

Next, we reformulate the MDP as an RMAB problem with simultaneous plays \cite{Whittle:1988JAP}.
The RMAB problem remains intractable in general and was shown to be PSPACE hard in \cite{Papadimitriou&Tsitsiklis:99INF}, which is in sharp contrast to the original (rested) multi-armed bandit (MAB) problem solved by the Gittin's index policy in \cite{Gittins&Jones:1979B}.
 Here we consider the celebrated Whittle's index policy that has been shown to be optimal in some special cases  \cite{Whittle:1988JAP,Liu&Qing:2010TIT}. To this end, we first establish the indexability of the formulated RMAB problem.  We then show that, for the deadline scheduling problem, in particular, the pre-determined deadline and workload at the time of arrival simplify the computation of the Whittle's index.  For the case with constant processing cost, we derive the Whittle's indexes in closed-form, which generalizes the result of \cite{Graczova&Jacko:14OR}.

When the number of processors is finite, we first provide examples that the Whittle's index policy is not optimal for the deadline scheduling problem. We show, however, that the gap-to-optimality for the Whittle's index policy is bounded by the conditional value at risk (CVaR) \cite{Rockafellar&Uryasev:2000JR} of the number of arrivals per unit time, which allows us to examine the performance loss as a function of arrival rate and number of available processors.

A major result of this paper is to characterize the asymptotic optimality of the Whittle's index policy when the number of processors increases with the job arrival rate.  In particular, we show that the gap-to-optimality  goes to zero in the light traffic case, indicating a specific regime in which Whittle's index policy is asymptotically optimal.

\subsection{Related Work}
The classical deadline scheduling problem is first considered by Liu and Layland \cite{Liu&Layland:73ACM} in a deterministic setting.  For the single processor case, the results are quite complete. When all jobs can be finished on time, simple index algorithms (with linear complexity) such as the earliest deadline first (EDF) \cite{Liu&Layland:73ACM,Dertouzos:74IFIPC} and the least laxity first (LLF) \cite{Mok:83Dissertation} achieve the same performance as the optimal off-line algorithm in the deterministic setting.

 There is also substantial literature on the
deadline scheduling problem with multiple processors
(for a survey, see \cite{DB11}). It is shown in \cite{DerTouzos&Mok:89TSE}
that optimal online scheduling policies do not exist in general for the worst case performance measure.

The literature on deadline scheduling in the stochastic settings is less extensive. For the single processor case, Panwar, Towsley, and Wolf in \cite{Panwar&Towsley&Wolf:88JACM} and \cite{Towsley&Panwar:90RT} made an early contribution in establishing  the optimality of EDF in minimizing the unfinished work when jobs are non-preemptive.  The performance of EDF is quantified in the heavy traffic regime using a diffusion model in \cite{Lehoczky:96RTSS,Doytchinov&Lehoczky&Shreve:01AAP}, and \cite{Kruk&Etal:11AAP}.

The multiprocessor stochastic deadline scheduling problem is less understood, primarily because the stochastic dynamic
programming for such problems are intractable to solve in
practice. A particularly relevant class of applications
is scheduling in wireless transmissions and routing in networks \cite{BT97,Jaramillo&Srikant&Ying:2011SAC,Raghunathan&Etal:2008INFOCOM}, and \cite{Singh&Kumar:15CDC} where job (packet) arrival
is stochastic, and packets sometimes have deadlines for
delivery. {Another class of applications
is in the scheduling of (deadline-constrained) electric vehicle charging with stochastic charging costs \cite{Xu&Pan:2012CDC,Jia15,XPT16}.}
The work closest to ours are in \cite{BT97,Jaramillo&Srikant&Ying:2011SAC,Raghunathan&Etal:2008INFOCOM}, and \cite{Singh&Kumar:15CDC} where the authors considered special instances of the deadline scheduling problem studied in this paper.  In the context of scheduling transmissions in wireless networks, the authors of \cite{BT97} analyzed the performance of the EDF policy for packets delivery in tree networks. Also related is the deadline scheduling in ad hoc networks \cite{Jaramillo&Srikant&Ying:2011SAC} where an iterative algorithm was proposed to schedule packets over random channels, and the algorithm was proved to be optimal. Random arrivals of jobs (packets) were considered in \cite{Raghunathan&Etal:2008INFOCOM} where the authors formulated the problem as an RMAB problem and analyzed the indexability.  Whittle's index policy was applied, but the performance of Whittle's index policy was not analyzed.  The model considered in \cite{Raghunathan&Etal:2008INFOCOM} is also more restrictive than the model studied in this paper. The work of \cite{Singh&Kumar:15CDC} considers the problem of scheduling multihop wireless networks for packets with deadlines where the authors developed decentralized scheduling policies. The constraint on bandwidth in \cite{Singh&Kumar:15CDC} is an average constraint whereas the problem treated in this paper is a strict deterministic constraint.

A recent related work in the operation research literature is \cite{Graczova&Jacko:14OR} where the authors considered the RMAB formulation of the deadline scheduling in knapsack problems.  The authors established the indexability of the RMAB problem and a closed-form of Whittle's index. There are several important differences, however, between the model considered in \cite{Graczova&Jacko:14OR} and the one in this paper.  First, the job arrivals are simultaneous in \cite{Graczova&Jacko:14OR} and stochastic in this paper.  Second, the processing cost/reward  is constant in \cite{Graczova&Jacko:14OR} and random in our model. Our paper also establishes the asymptotic performance of Whittle's index policy whereas \cite{Graczova&Jacko:14OR} addressed the indexability and developed an iterative algorithm to compute the Whittle's index.

There is extensive literature on the RMAB problem.  See, e.g., \cite{Ruiz:book,Gittins&Etal:2011book}.  In his seminal work \cite{Whittle:1988JAP}, Whittle introduced an index policy (the Whittle's index policy) for the subclass of indexable RMAB problems. Although in general suboptimal in the finite arm regime  except for some special cases \cite{Liu&Zhao:08ITsub}, Whittle's  index policy was shown by Weber and Weiss in \cite{Weber&Weiss:JAP1990} to be asymptotically optimal under some conditions when the number of arms and the number of simultaneous activations grow proportionally to infinity. The optimality conditions, however, are difficult to check.  We should also point out that the asymptotic optimality results established in this paper is different from that formulated in \cite{Whittle:1988JAP} and \cite{Weber&Weiss:JAP1990}.

%
%
%

%% file: formulation_v8.tex
In this section, we introduce the stochastic deadline scheduling problem as a constrained MDP followed by an RMAB formulation.

 \subsection{Nominal Model Assumptions}\label{sec:IIA}
We begin with a set of nominal assumptions in setting up the MDP formulation:
\bitem
\item [A1.]The time is slotted, indexed by $t$.
\item[A2.] There are $M$ processors available at all times. In each time slot, a processor can only work on one job, and each job can receive service from only one processor at any given time.  A processor can be switched from one job to another without incurring  switching cost.
\item[A3.] If a processor works on a job in time slot $t$, it receives a unit payment and incurs a time-varying cost $c[t]$. Here we assume that $c[t]$ is an exogenous stationary Markov process with a transition probability matrix ${P=[P_{i,j}}]$. 
\item[A4.] If a job is not completed by its deadline, a penalty defined by a convex function  of the amount of unfinished job is imposed on the scheduler by the deadline.
    Let $F(B)$ be the convex penalty function with $B$ denoting the amount of the unfinished job with ${F(0)=0}$.

    \item[A5.] There is a queue with $N$ positions. The jobs arrived at different positions are statistically independent and identically distributed (i.i.d.). %

\item [A6.]   A job arriving at the $i$th position of the queue at the beginning of time slot $t$ reveals $B_i$ (the total amount of work to be completed) and $T_i$ (the deadline for completion). At the end of time slot $t+T_i$, the job is removed from the queue, regardless whether the job is completed. When the $i$th position is available, with probability $Q(T,B)$ a new job with deadline $T$ and workload $B$ arrives.  With probability $Q(0,0)$, the position remains empty.

\eitem

Some comments and clarifications on these assumptions are in order. Assumptions A1 and A2 are standard.
 A3 assumes that the marginal price of service---the marginal payment to the service provider---is the same for all jobs. The marginal processing cost $c[t]$ is uniform for jobs processed at the same $t$.  Several generalizations of A3 are possible.  In particular, by including the initial lead time in the state of a job, our model can accommodate the so-called service differentiated deadline scheduling problem \cite{bitar2014deadline}, where jobs with different deadlines face different marginal prices.  Another generalization is that the marginal price (or the cost) of service depends on the position of the queue.  This, for instance, can model prioritized services.

Assumption A4 indicates that the deadline is soft, but it can be hardened by setting the non-completion penalty much higher than processing cost.
In this setting, it is always optimal ({\it i.e.}, reward maximizing) to finish as many jobs as possible, regardless of the processing cost.

The i.i.d. arrival assumption in A5 is limiting but necessary for index policies.  This is also consistent with the standard Poisson arrival case when the arrived job is randomly assigned to a position in the queue.  A5 and A6 imply that when a job arrives at a position that is occupied by an unfinished job, the newly arrived job is dropped, which seems unreasonable since the job could have been reassigned to an open position (if it exists). However, asymptotically when ${N \rightarrow \infty}$, there is no loss of performance by imposing these assumptions. In Section \ref{sec:validationOfA6}, we numerically compare the two scenarios with i.i.d. arrivals following A5-A6 and the conventional Poisson arrival. Numerical results show that the performance of different algorithms
under A6 converges to its counterpart under Poisson arrival as the number of available positions increases.
\subsection{Stochastic Deadline Scheduling as a Constrained MDP}
Next, we formulate the constrained MDP by defining states, actions of the scheduler, state evolution, rewards, constraints, and decision policies.

\subsubsection{State Space}
\begin{figure}[h]
\centering
\psfrag{time}{{\scriptsize time}}
\psfrag{t}{{\scriptsize$t$}}
\psfrag{car1}{{\scriptsize Job $J_i$}}
\psfrag{rn}{{\scriptsize $r_i$}}
\psfrag{Tn}{{ \scriptsize $T_i[t]$}}
\psfrag{Bn}{{\scriptsize$B_i[t]$}}
\psfrag{xn}{{\scriptsize$L_i[t]$}}
\psfrag{dn}{{\scriptsize $d_i$}}
  \includegraphics[width=.45\textwidth]{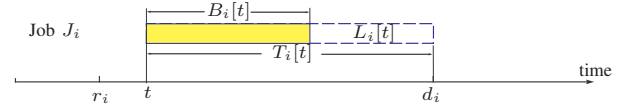}
  \caption{An illustration of job $i$'s state.  $r_i$ is the arrival time of a job at position $i$, $d_i$ is its deadline for completion, $B_i[t]$ is the workload to be completed by $d_i$, $T_i[t]$ is the job's lead time to deadline, and $L_i[t]\triangleq T_i[t]-B_i[t]$ is the job's laxity.}
  \label{fig:carExample}
\end{figure}

Consider first the state of  the $i$th position in the queue.  Let ${T_i[t]\triangleq d_i-t}$ be the lead time to deadline $d_i$, $B_i[t]$ be the remaining job length, and $L_i[t]\triangleq T_i[t]-B_i[t]$ be the laxity of job $i$, as illustrated in Figure~\ref{fig:carExample}.

The state of  the $i$th position in the queue is defined as
\[
S_i[t] \defeq \left\{\begin{array}{ll}
(0,0) & \mbox{if no job is at the $i$th position,}\\
(T_i[t],B_i[t]) & \mbox{otherwise.}\\
\end{array}
\right.
\]
The processing cost $c[t]$ evolves according to an exogenous finite state Markov chain with a transition probability matrix ${P=[P_{j,k}}]$. This Markovian assumption is practical to study stochastic prices, \eg \cite{Kakade&Kearns:2005ICCLT}, and simplifies both the model and the computation of the policies.

The state of the MDP is defined by the queue states and the processing cost $c[t]$ as
${S[t]\defeq(c[t], S_1[t],\cdots, S_N[t])\in\mathcal{S}}$, where ${\mathcal{S}}$ is the state space.

\subsubsection{Action}
The action of the scheduler in slot $t$ is defined by the binary vector $\abf[t]=(a_1[t], \cdots, a_N[t]) \in \{0,1\}^N$.
 When
 $a_i[t]=1$, a processor is assigned to work on the job at position $i$, and the position is called {\em active}.
 When $a_i[t]=0$, position $i$ is {\em passive}, \ie no processor is assigned.   For notational convenience, sometimes we allow a position without a job to be activated, in which case the assigned processor  receives no reward and incurs no cost.

 \subsubsection{State Evolution}\label{sec:stateEvolution}
The evolution of the processing cost is according to the transition matrix $P$ and
 independent of the actions taken by the scheduler. The evolution of the job state $S_i[t]$  depends on the scheduling action $a_i[t]$:
\beq
S_i[t+1] = \left\{\begin{array}{ll}
(T_i[t]-1,(B_i[t]-a_i[t])^+)& \mbox{if }T_i[t]>1, \\
(T,B)~\mbox{with prob. $Q(T,B)$} & \mbox{if }T_i[t] \le 1, \\
\end{array}
\right.
\label{eqn:stateTransitionActive}
\eeq
where ${b^+=\max(b,0)}$.
Note that when ${T_i[t]=1}$, the deadline is due at the end of the current time slot
 and the job in position $i$ will be removed. 

\subsubsection{Reward}
 For each job, the scheduler obtains one unit of reward if the job is processed for one time slot. When ${T_i[t]=1}$, job $i$ will reach its deadline by the end of the current time slot, and the scheduler will incur a penalty if the job is unfinished.
The reward collected from job $i$ at time $t$ is given by
\beq\hspace{-1em}
\begin{array}{l}
R_{a_i[t]}(S_i[t],c[t])\\
~~~~=\left\{\begin{array}{ll}
(1-c[t])a_i[t] & \mbox{if }B_i[t]>0, T_i[t]>1,\\
(1-c[t])a_i[t] \\
-F(B_i[t]-a_i[t]) & \mbox{if }B_i[t]>0, T_i[t]=1,\\
 0 & \mbox{otherwise.}
 \end{array}
 \right.
 \end{array}
 \label{eqn:reward}
 \eeq

\subsubsection{Objective}
 Given the initial system state ${S[0]=s}$ and a policy $\pi$ that maps
 each system state $S[t]$ to an action vector $\textbf{a}[t]$,
 the expected discounted system reward is defined by
 \beq \label{eq:V}
G^N_\pi(s) \, \defeq \,\mathbb{E}_\pi \left(\sum_{t=0}^\infty \sum_{i=1}^N \beta^t R_{a_i[t]}(S_i[t],c[t]) \bigg|  S[0]=s \right),
 \eeq
 where $\mathbb{E}_\pi$ is the conditional expectation over the randomness in costs and job arrivals under a given scheduling policy $\pi$ and $0<\beta<1$ is the discount factor. 

 \subsubsection{Constrained MDP and Optimal Policies}
We impose a constraint on the maximum number of processors that can be activated simultaneously, \ie $\sum_{i=1}^N a_i[t] \le M$.  This constraint represents the processing capacity of the service provider. For the EV charging application, this assumption translates directly to the physical power limit imposed on the charging facility.
  Thus, the deadline scheduling problem can then be formulated as a constrained MDP.
\beq
\label{eqn:originPro}
G^N(s)=\sup_{\{\pi: \sum_{i=1}^N a^{\pi}_i[t] \le M, \, \forall t\}}  G^N_\pi(s),
\eeq
where $a^{\pi}_i[t]$ is the action sequence generated by policy $\pi$ for position $i$.  A policy $\pi^*$ is optimal if $G^N_{\pi^*}(s)=G^N(s)$.
Without loss of optimality, we will restrict our attention to  stationary policies \cite{Altman:2004book}.

\subsection{A Restless Multi-armed Bandit  Problem}
Unfortunately, the MDP formulation does not result in a scalable optimal scheduling policy because the state space grows exponentially with $N$.
We, therefore, seek to obtain an effective {\it index policy} \cite{Gittins:79JRSS} that scales linearly with $N$. We identify each position in the queue as an arm and formulate (\ref{eqn:originPro}) as an RMAB problem. To this end,  ``playing'' an arm is equivalent to assigning a processor to process the job (if there is one) at a position in the queue.
  The resulting multi-armed bandit problem is restless because the state of position $i$---the $i$th arm---evolves regardless whether arm $i$ is active or passive.

A complication of casting  (\ref{eqn:originPro}) as an RMAB problem comes from the inequality constraint on the maximum number of simultaneously activated positions,
as the standard RMAB formulation imposes an equality constraint on the number of arms that can be activated.
This can be circumvented by introducing  $M$ {\em dummy arms} and requiring that exactly $M$ arms must be activated in each time slot.  Specifically, each dummy arm $i$ always accrues zero rewards, and its state stays at ${S_i=(0,0)}$.   The reformulated RMAB problem has ${N+M}$ arms.  We let ${\{1,\cdots,N\}}$ be the set of \emph{regular arms} that generate reward (penalty) and ${\{N+1,\cdots,N+M\}}$ be the set of
dummy arms.

 We define the extended state of each arm as ${\tilde{S}_i[t]\triangleq(S_i[t],c[t])}$ and denote the extended state space as ${\Smsc_i\triangleq\mathcal{S}_i\times \mathcal{S}_c}$. The state transition of each arm and the associated reward inherit from (\ref{eqn:stateTransitionActive}-\ref{eqn:reward}) of the original MDP.
 We have the following RMAB problem that is equivalent to the original MDP (\ref{eqn:originPro}):
\begin{equation}
\label{eqn:MAB}
  \begin{array}{ll}
\sup_{\pi} & \mathbb{E}_\pi\left\{\sum_{t=0}^{\infty}\sum_{i=1}^{N+M}\beta^{t}R_{a_i[t]}(\tilde{S}_i[t])\mid \tilde{S}_i[0]\right\}\\[5pt]
\mbox{s.t.}& \sum_{i=1}^{N+M}a_i[t]= M, \quad {\forall } \, t.
  \end{array}
  \end{equation}
In (\ref{eqn:MAB}), the arms are coupled by the processing cost.
With the addition of dummy arms, the inequality constraint on the maximum number of activated arms in the original MDP problem is transformed to the equality constraint
in (\ref{eqn:MAB}).

%% file: whittleIndex_v7.tex
To tackle the deadline scheduling problem as an RMAB, we first establish the indexability of the RMAB
and then formally define the Whittle's index policy in this section.

\subsection{Indexability}
Following \cite{Whittle:1988JAP}, we consider the {\em $\nu$-subsidized} single arm reward maximization problem that
seeks for a policy $\pi$ to activate/deactivate a single arm to maximize the discounted accumulative reward:
\begin{equation}
V^\nu_i(s)=\sup_{\pi}\mathbb{E}_\pi\left(\sum_{t=0}^{\infty}\beta^{t}R^\nu_{a_i[t]}(\tilde{S}_i[t])\bigg| \tilde{S}_i[0]=s \right),
\label{eq:nu-subsidy}
  \end{equation}
where the subsidized reward is given by
\[
R^\nu_{a_i[t]}(\tilde{S}_i[t])=R_{a_i[t]}(\tilde{S}_i[t])+\nu \mathds{1}(a_i[t]=0).
\]
Here  $R_{a_i[t]}(\cdot)$  is defined in (\ref{eqn:reward}),
and $\mathds{1}(\cdot)$ is the indicator function.
In the {\em $\nu$-subsidized} problem, the  scheduler receives a subsidy $\nu$
whenever an arm is passive.

Let $\mathcal{L}_{a}$ be an operator on $V_i^{\nu}$ defined by
\[
(\mathcal{L}_{a}V_i^\nu)(s)
\triangleq\mathbb{E}\left(V_i^{\nu}(\tilde{S}_i[t+1])\bigg| \tilde{S}_i[t]=s,a_i[t]=a\right).
\]
The maximum discounted reward $V^\nu_i(\cdot)$ in (\ref{eq:nu-subsidy}) is determined by the Bellman equation
\beq\label{eqn:BellmanNuSubsidy}
V^\nu_i(s)=\max_{a \in \{0,1\}}\Bigg\{R^\nu_a(s)+\beta(\mathcal{L}_{a}V^\nu_i)(s)\Bigg\}.
\eeq

Let  $\Smsc_i(\nu)$ be  the set of states under which it is optimal to take the passive action in the $\nu$-subsidy problem. The {\em indexability} of the RMAB is defined by the monotonicity of $\Smsc_i(\nu)$ as the subsidy level $\nu$ increases.

\begin{Def}[Indexability \cite{Whittle:1988JAP}]\label{def:index}
  Arm $i$ is indexable if the set $\Smsc_i(\nu)$ increases monotonically from $\emptyset$ to $\Smsc_i$ as $\nu$ increases from $-\infty$ to $+\infty$. The MAB problem is indexable if all arms are indexable.
\end{Def}

We establish the indexability for the stochastic deadline scheduling problem.

\begin{thm}[Indexability]\label{thm:indexability}${}$Each arm is indexable, and the RMAB problem (\ref{eqn:MAB}) is indexable.
\end{thm}
The indexability of the bi-dimension state model without arrival is proved in \cite{Graczova&Jacko:14OR} based on the partial conservation law principle. We provide an elementary proof in Appendix \ref{proof:indexability} that also includes the random arrivals of jobs.

\subsection{Whittle's Index Policy}
 The following definition of Whittle's index is based on Definition \ref{def:index}.

\begin{Def}[Whittle's index \cite{Whittle:1988JAP}]\label{def:W}
If arm $i$ is indexable, its Whittle's index $\nu_i(s)$ of state $s$ is the infimum of the subsidy $\nu$ under which the passive action is optimal at state $s$,
{\it i.e.},
\[
\begin{array}{l}
\nu_i(s)\triangleq\inf_\nu\{\nu:R_0(s)+\nu+\beta(\mathcal{L}_{0}V_i^\nu)(s)\\[3pt]
\mathrel{\phantom{\nu_i(\tilde{s})\triangleq\inf_\nu\{\nu:}}\ge R_1(s)+\beta(\mathcal{L}_{1}V_i^\nu)(s)\}.
\end{array}
\]
\end{Def}
If arm $i$ is indexable, in a $\nu$-subsidized problem with ${\nu<\nu_i({s})}$  it is optimal to activate arm $i$.
Likewise, if $\nu\ge\nu_i(s)$  it is optimal to  deactivate arm $i$.

To compute the Whittle's index for arm $i$, we solve a parametric program where the subsidy $\nu$ appears in the constraints.
 \[
 \begin{array}{ll}
 \min_{u_i(s)}& \sum_{s\in\mathcal{S}}p(s)u_i(s)\\[2pt]
 \mbox{s.t.}& u_i(s)\ge R_{1}(s)+\beta \sum_{s'\in\mathcal{S}}P_{s,s'}^1u_i(s'), \quad {\forall } s,\\[2pt]
& u_i(s)\ge R_{0}(s)+\nu+\beta \sum_{s'\in\mathcal{S}}P_{s,s'}^0u_i(s'),\quad {\forall } s,\\[2pt]
&~~~~~~~~~~~~~
 \end{array}
 \]
where ${s=(T,B,c)}$ is the extended state of arm $i$, $p(s)$ the initial state probability, and $P_{s,s'}^a$ the transition probability from $s$ to $s'$ given action $a$. For a particular value of $\nu$, the optimal solution $u^*_i(s)$ equals the value function $V^\nu_i(s)$,
and one of the two active constraints gives the optimal action. We solve this parametric program to find the break point of $\nu$ where the optimal action changes.
The simplex method can be used to solve this parametric program \cite{Dantzig:book1998}.

The special structure of the deadline problem allows us to obtain a closed-form solution when the processing cost is time-invariant.

\begin{thm}\label{thm:closed-form}
If $c[t]=c_0$ for all $t$, the Whittle's index of a regular arm $i\in\{1,\cdots,N\}$ is given by
\begin{equation}\label{eqn:closedForm}
\begin{array}{l}
\nu_i(T,B,c_0)\\
~~~~=\left\{
\begin{array}{ll}
0 &\mbox{if~}\, B=0,\\[3pt]
1-c_0 &\mbox{if~}\, 1\le B\le T-1,\\[3pt]
\beta^{T-1}F(B-T+1)\\
-\beta^{T-1}F(B-T)\\
+1-c_0&\mbox{if~}\, T\le B.
\end{array}
\right.
\end{array}
\end{equation}
The Whittle's index of a dummy arm is zero, \ie
\[
\nu_i(0,0,c_0)=0, \quad i\in\{N+1,\cdots,N+M\}.
\]
\end{thm}

 The proof of Theorem \ref{thm:closed-form} can be found in Appendix \ref{proof:closed-form}.
In (\ref{eqn:closedForm}), when it is feasible to finish job $i$'s request ({\it i.e.} its lead time is no less than its remaining processing time),
job $i$'s Whittle's index is simply the (per-unit) processing profit $1-c_0$. When a non-completion penalty is inevitable, the index takes into account both
the processing profit and the non-completion penalty. We note that the Whittle's index gives higher priority to jobs with less laxity.

We are now ready to define the Whittle's index policy based on Definition \ref{def:W}.

\begin{Def}[Whittle's index policy \cite{Whittle:1988JAP}]\label{def:WhittleIndexPolicy}
For the RMAB problem defined in (\ref{eqn:MAB}), the Whittle's index policy sorts all arms by their Whittle's indices  in a descending order and activates the first $M$ arms.
\end{Def}

Since the states of jobs and processing cost are finite,
the Whittle's indices can be computed off-line. In real-time
scheduling, at the beginning of each time slot, the scheduler
looks up the indices for all existing jobs based on the current system state and processes the ones with
highest indices. When there is a tie, the scheduler breaks the
tie randomly with a uniform distribution.

We note that the Whittle's index policy does not distinguish jobs with positive laxity,
which leaves some room for improvement. In Section \ref{sec:op4}, we apply the
Least Laxity and Longer Processing Time (LLLP) principle (originally proposed in \cite{Xu&Pan:2012CDC}) to
improve the Whittle's index policy.

%% file: finiteArm.tex
In this section, we examine the performance of Whittle's index policy for the stochastic deadline scheduling problem when the number of servers ($M$) is finite.  We show that when ${M<N}$, there does not exist an optimal index policy in general, hence Whittle's index policy is not optimal.   We further derive an upper bound on the gap-to-optimality on the performance of the Whittle's index policy.  This result provides the essential ingredient for establishing asymptotic optimality of the Whittle's index policy in Section \ref{sec:V}.

%
%
%
%
%

\subsection{Performance in the Finite Processor Cases}\label{sec:op1}
In general, Whittle's index policy is not optimal except in some special cases \cite{Liu&Qing:2010TIT}.  For the deadline scheduling problem, the same conclusion holds. We show in fact that no optimal index policy exists.

\begin{prop}\label{prop:optimal}
 When ${M=N}$, the Whittle's index policy is optimal. When ${M<N}$, an optimal index policy for the RMAB problem formulated in (\ref{eqn:MAB}) may not exist in general.
\end{prop}
\begin{proof} The fact that Whittle's index policy is optimal when ${M=N}$ is intuitive. A formal proof can be found in Appendix \ref{proof:optimalityOfWhittleIndexPolicy}. To show that an optimal index policy does not exist in general, it suffices to construct a counter example that no index policy can be optimal.

Set the capacity of the queue to be ${N=3}$, the number of processors ${M=1}$, the discounted factor ${\beta=0.4}$, the penalty function ${F(B)=B^2}$, and the processing cost ${c[t]=1}$. Assume the arrival is busy $(Q(0,0)=0)$ and the initial laxity is zero ($T=B$ at arrival). For this small scale MDP, a linear programming formulation is used to solve for the optimal policy \cite{Ross:book2014}.

 Consider two different states,
 \[
 \begin{array}{ll}
 {s=((1,1),(2,2),(2,2))},\\
 {s'=((1,1),(1,1),(2,2))},
 \end{array}
 \]
where ${s=((T_1,B_1),(T_2,B_2),(T_3,B_3))\in\mathcal{S}}$ is the state of the system including the states of each arm. The constant processing cost is omitted in the state.

 For state $s$, the optimal action is to process job $(2,2)$. In this case, the job $(2,2)$ is preferred to $(1,1)$. Processing $(2,2)$ will cause an immediate  penalty of $1$, and the state will change to ${((T,B),(1,1),(1,2))}$, where $(T,B)$ is a new arrival. In next stage, a penalty of $2$ from the last two jobs will happen. If some policy processes $(1,1)$ alternately given state $s$, there will be no penalty in the first stage, and the state will change to ${((T,B),(1,2),(1,2))}$. The last two jobs will at least incur a penalty of $5$.

 For state $s'$, the optimal action is to process the job $(1,1)$. The job $(1,1)$ is preferred to $(2,2)$ in this case. Processing $(1,1)$ will cause an instant penalty of $1$, and the state will change to ${((T,B),(T',B'),(1,2))}$, where $(T,B)$ and $(T',B')$ are new arrivals. If some policy processes $(2,2)$ alternately given state $s'$, there will be an instant penalty of $2$ from the first two jobs in the first stage and the state will change to ${((T,B),(T',B'),(1,1))}$. In this case, a penalty of $1$ can be saved in the second stage by processing $(2,2)$ in the first one. However, due to the discount factor, it is more profitable to process $(1,1)$.

 An {\it index policy} assigns each job an index (that depends only on the job's current state) and processes the jobs with the highest indices \cite{Gittins:79JRSS}. Therefore, for any ``index'' policy, the indices of job $(1,1)$ and $(2,2)$ are fixed, and the preference of these two jobs should remain the same in these two cases, which is violated by the result here. This counter example shows that no ``index'' policy that is optimal in general.
\end{proof}

Note that, the Whittle's index policy is an example of index policies, and thus is sub-optimal in general. However, with particular combinations of parameters, optimal index policies may exist.

\subsection{An Upper Bound of the Gap-to-Optimality}\label{sec:op3}
In the following lemma, we first establish a result that applies quite generally to the case for a  finite queue size $N$ and a finite number of processors $M$.

\begin{lem}\label{lem:gap}
Let $G^N(s)$ be the optimal value function defined in (\ref{eqn:originPro}) and $G^N_{\scriptsize\mbox{W}}(s)$ be the value function achieved by the Whittle's index policy, respectively.  We have
  \beq\label{eqn:bound}
  \begin{array}{l}
  G^N(s)-G^N_{\scriptsize\mbox{W}}(s)\\
  ~~~~\le \frac{C}{1-\beta}\mathbb{E}[I^N[t]|I^N[t]>M]\Pr(I^N[t]>M),
  \end{array}
  \eeq
  where $I^N[t]$ is the number of jobs admitted in the queue with $N$ positions within time $[t-\bar{T}+1,t]$, $\bar{T}$ is the maximum lead time of jobs, and $C$ is a constant determined by the processing cost and the penalty of non-completion.
  \end{lem}

The proof can be found in Appendix \ref{proof:asymptotic}. The gap-to-optimality is bounded by the tail expectation of the jobs admitted to the system. Note that, the conditional expectation on the right-hand side (RHS) of (\ref{eqn:bound}) is connected to the conditional value at risk (CVaR) \cite{Rockafellar&Uryasev:2000JR}, which measures the expected losses at a certain risk level and is extremely important in the risk management.

\subsection{Least Laxity and Longer Processing Time (LLLP) Principle}\label{sec:op4}
\input scheduling_v4.tex

%% file: scheduling_v4.tex
In this subsection, we will apply the Less Laxity and Longer remaining Processing time (LLLP) principle (originally proposed in \cite{Xu&Pan:2012CDC}) to improve the Whittle's index policy. As a priority rule for stochastic deadline scheduling, the LLLP principle is defined as follows.

\begin{Def}[LLLP Order \cite{Xu&Pan:2012CDC}] \label{def:LLLP}
Consider jobs $i$ and $j$ at time $t$. We say $j$ dominates $i$ ($j\succeq i$)
if $j$ has less laxity and longer remaining processing time than those of $i$, \ie  $L_j[t]\le L_i[t]$ and ${B_j[t] \ge B_i[t]}$, with at least one of the inequalities strictly holds.
\end{Def}

LLLP defines a partial order over the jobs' states such that the job with less laxity and longer remaining job length should be given priority. Compared to LLF, LLLP takes into account both the laxity and the remaining workload, whereas LLF considers laxity only.

An LLLP interchange enhancement policy is proposed in \cite{Xu&Pan:2012CDC}. Specifically, it is shown that applying the LLLP
interchange on a policy $\pi$ leads to a policy that performs
no worse than that of $\pi$.
 Numerical experiments shown in \cite{Xu&Pan:2012CDC} demonstrates that LLLP enhancement often performs
significantly better than the policy to which LLLP is applied.
  This insight leads us to consider an LLLP enhancement on
the Whittle's index policy in the context of RMAB approach
to stochastic deadline scheduling.

Denote the set of arms by ${\mathcal{N}=\{1, \cdots,N+M\}}$. Consider a policy $\pi$  that activates arms (jobs) in $\Xmsc$ and deactivates those in ${\Xmsc^c=\mathcal{N}\setminus \Xmsc}$ at system state $S[t]$. We say that a policy $\pi$ follows the LLLP principle if there does not exist a pair of jobs ${(i,j)}$ such that ${i\in\Xmsc}$, ${j\in\Xmsc^c}$, and ${j\succeq i}$. We propose a Whittle's index based algorithm that activates arms with highest indices without violating the LLLP principle.

\begin{figure}[ht]

\psfrag{1}{{ $1$}}
\psfrag{2}{{ $2$}}
\psfrag{3}{{ $3$}}
\psfrag{4}{{ $4$}}
\psfrag{5}{{ $5$}}
\psfrag{6}{{ $6$}}
\psfrag{7}{{ $7$}}
\psfrag{8}{{ $8$}}
\psfrag{9}{{ $9$}}
\psfrag{10}{{ $10$}}
  \includegraphics[width=.48\textwidth]{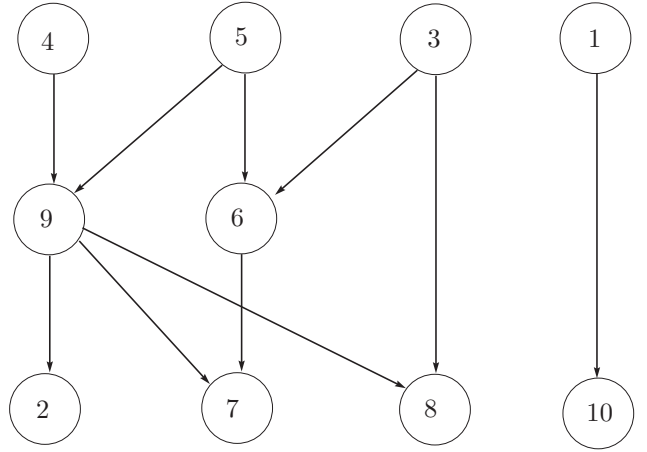}
\caption{A directed acyclic graph indicating the LLLP order.}
\label{fig:toposort}
\end{figure}

As shown in Figure \ref{fig:toposort}, the LLLP order defines a directed acyclic graph (DAG) ${\mathcal{G}=\{\mathcal{N},\mathcal{E}\}}$ of all arms, where $\mathcal{N}$ represents the arm set and $\mathcal{E}$ is the edge set. An edge from $i$ to $j$ indicates that job $i$ dominates job $j$ in the sense of LLLP order. A topological sorting is a linear ordering of the vertices so that for each directed edge ${(i,j)\in\mathcal{E}}$, $i$ comes before $j$ in the ordering.

Typically, topological sorting of a DAG is not unique. We employ a stable topological sorting to guarantee that the result ordering preserves the order of Whittle's index of arms whenever it is possible. In the proposed algorithm, we employ a depth-first search with linear complexity in the number of vertices and edges \cite{Thomas&Etal:2001book}. In Figure \ref{fig:toposort}, arms are pre-ordered descendingly according to their Whittle's indices, and the LLLP ordering is indicated by the directed edges\footnote{In the proposed algorithm, if an arm has a state of $(0,0)$ (either no job or a dummy arm), it dominates no arm and no arm dominates it.}. The stable topological sorting gives an order of $\{1,3,4,5,6,9,2,7,8,10\}$. The LLLP enhanced Whittle's index policy is formulated in Algorithm \ref{alg:WhittleLLLP}.

%
%
%
%

\begin{algorithm}
\caption{Whittle index with LLLP interchange}
\label{alg:WhittleLLLP}
\begin{algorithmic}
\STATE $1$. Calculate the Whittle's indices of all arms and sort them in a descending order.
\STATE $2$. Generate a DAG according to the LLLP ordering.
\STATE $3$. Carry out a stable topological sort.
\STATE $4$. Activate the $M$ arms with the highest priority.
\end{algorithmic}
\end{algorithm}

%% file: asymptotic_v8.tex
In this section, we establish the asymptotic optimality of the Whittle's index policy when the job arrival rate $\mu$ and the number of servers $M$ increase to infinity simultaneously while the system stays stable.

We first consider the case when the aggregated arrival of jobs follows a Poisson distribution. Let $I[t]$ be the total number of jobs arrived at the system within $[t-\bar{T}+1, t]$, recalling that $\bar{T}$ is the maximum lead time of jobs.   Note that $I[t]$ is Poisson distributed.

When the queue at the service center is finite with $N$ positions, we assume that each position receives equally likely $1/N$th of the traffic\footnote{The thinning property of Poisson justifies A6.}.   Because a newly arrived job may be rejected when the assigned position is occupied (A5), the total number of jobs $I^N[t]$ admitted to the system in slot $t$ satisfies $I^N[t]\le I[t]$. Moreover, we have the following corollary.
\begin{corollary}\label{col:infiniteN}
As $N \rightarrow \infty$, $I^N[t]\rightarrow I[t]$ in distribution. 
\end{corollary}
The proof can be found in Appendix \ref{proof:infiniteN}.

 Define
\[
G(s)-G_{\scriptsize\mbox{W}}(s)\triangleq \limsup_{N\rightarrow\infty}[G^N(s)-G^N_{\scriptsize\mbox{W}}(s)].
\]
Then, by Lemma \ref{lem:gap} and Corollary \ref{col:infiniteN},
    \begin{equation}\label{eqn:infiniteN}
G(s)-G_{\scriptsize\mbox{W}}(s)\le\frac{C}{1-\beta}\mathbb{E}[I[t]\mathds{1}(I[t]>M)].
\eeq

Equation (\ref{eqn:infiniteN}) characterizes the performance gap for the Whittle's index policy for
 the asymptotic regime as $N$ increases while the arrival process and number of processors stay constant. Now, we check the performance of the Whittle's index policy when the number of processors $M$ increases and the mean of the arrival process $I[t]$  also grows as a function.

\begin{thm}\label{thm:Poisson}  Suppose that the aggregated arrival $I[t]$ is Poisson with mean $\mu$.  The Whittle's index policy is asymptotically optimal as $M \rightarrow \infty$ if $\mu < M/e$.  In particular,
\beq\label{eqn:Poisson}
G(s)-G_{\scriptsize\mbox{W}}(s) = \mathcal{O}(\frac{\mu e^{-\mu}}{\sqrt{M}}).
\eeq
\end{thm}

The proof of Theorem \ref{thm:Poisson} can be found in Appendix \ref{proof:Poisson}.  Besides showing that the Whittle's index is asymptotically optimal, Theorem \ref{thm:Poisson} also indicates that the gap-to-optimality decays sub-exponentially when $\mu$ grows with $M$ at the constant rate less than $1/e$.  When $\mu$ grows slower than $M$, the gap decays to zero but at a slower rate.

In general, suppose that we don't have the aggregated Poisson arrival, but  $I^N[t]$ converges in distribution to ${\bar{I}[t]\le I[t]}$ as $N\rightarrow \infty$.  If $I[t]$ with mean $\tilde{\mu}$ has a light tailed distribution, \ie there
exist  constants $a\ge1$ and $b\ge 0$ with
\begin{equation}\label{eqn:lightTail}
\Pr(I[t]\ge i)\le a\exp[-ib/\tilde{\mu}], \; \forall i\ge0,
\end{equation}
it can then be shown in \cite{Yu:17Dissertation} that,
\beq\label{eqn:exp}
G(s)-G_{\scriptsize\mbox{W}}(s)= \mathcal{O}[\exp(-\frac{Mb}{\tilde{\mu}})(Mb+\tilde{\mu})],
\eeq
as $M\rightarrow \infty$.

If $I[t]$ has a heavy tailed distribution with mean $\tilde{\mu}$, \ie
there exist constants $a>0$ and ${b>2}$
with
\beq\label{eqn:heavyTail}
\Pr(I[t]\ge i)\le a\tilde{\mu} /i^{b}, \forall i>0,
\eeq
it can then be shown in \cite{Yu:17Dissertation} that,
\beq\label{eqn:polynomial}
G(s)-G_{\scriptsize\mbox{W}}(s)=\mathcal{O} (\tilde{\mu}/M^{b-1}),
\eeq
as $M\rightarrow \infty$.

In both cases, the Whittle's index policy is asymptotically optimal if the arrival rate grows at the order of $o(M)$.

%% file: simulation_v8.tex
In this section, we present numerical results to compare  the performance of the Whittle's index policy  with other simple heuristic (index) policies,
 {\it i.e.}, EDF (earliest deadline first) \cite{Liu&Layland:73ACM}, LLF (least laxity first) \cite{Dertouzos:74IFIPC}, and Whittle's index policy with LLLP enhancement (cf. Algorithm~\ref{alg:WhittleLLLP}).

 If feasible, EDF processes $M$ jobs with the earliest deadlines, and  LLF processes $M$ jobs with the least laxity. Both algorithms break ties randomly.
Note that both policies will fully utilize the processing capacity and activate $M$ jobs as long as there are at least $M$ unfinished jobs in the system.
The Whittle's index policy, on the other hand,  ranks all arms by the Whittle's index and activates the first $M$ arms, and may
put some (regular) arms idle (deactivated) when the processing cost is high. The performance upper bound was obtained by replacing the strict capacity limit constraint by the constraint on the average \cite{Whittle:1988JAP}.

\subsection{Time-invariant Processing Cost}
We first considered a special case of problem (\ref{eqn:MAB})
  with a constant processing cost.
    Since the processing cost was time-invariant,
    it was optimal to fully utilize the capacity to process $M$ unfinished jobs.

\begin{figure}[ht]

\psfrag{Num of arm}{{\small $M/N$}}
\psfrag{Reward per arm}{{\small Total reward/N (\$)}}
\psfrag{UB}{{\scriptsize Upper Bound}}
\psfrag{EDF}{{\scriptsize EDF}}
\psfrag{LLF}{{\scriptsize LLF}}
\psfrag{Valley filling}{{\scriptsize Valley filling}}
\psfrag{Whittls index}{{\scriptsize Whittle's index}}
\psfrag{Whittle's index+LLLP}{{\scriptsize Whittle's+LLLP}}
  \includegraphics[width=.48\textwidth]{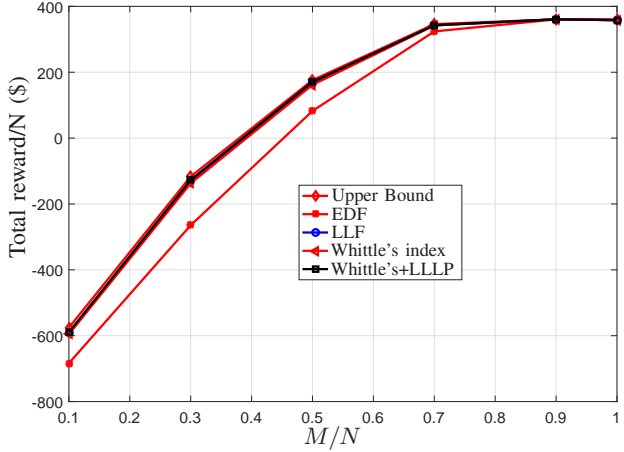}
\caption{Performance comparison with constant processing cost: ${c[t]=0.5}$, ${Q(0,0)=0.3}$, ${\bar{T}=12}$, ${\bar{B}=9}$, ${\beta=0.999}$, ${F(B)=0.2B^2}$, ${N=10}$.}
\label{fig:constantCostVaryRatioReward}
\end{figure}

In Figure~\ref{fig:constantCostVaryRatioReward}, we fixed the job arrival process and the length of the queue $N$ and varied the processing capacity $M$.
All policies except the EDF policy performed well and achieved an expected reward
close to the performance upper bound. When ${M/N=1}$, all jobs could be finished, and all policies achieved optimality.

\begin{figure}[ht]
\centering
\psfrag{Num of arm}{{\small $N$}}
\psfrag{Reward per arm}{{\small Total reward/N (\$)}}
\psfrag{UB}{{\scriptsize Upper Bound}}
\psfrag{EDF}{{\scriptsize EDF}}
\psfrag{LLF}{{\scriptsize LLF}}
\psfrag{Valley filling}{{\scriptsize Valley filling}}
\psfrag{Whittl's index}{{\scriptsize Whittle's index}}
\psfrag{LLLLLLLLLLLLP-Index}{{\scriptsize Whittle's+LLLP}}
  \includegraphics[width=.48\textwidth]{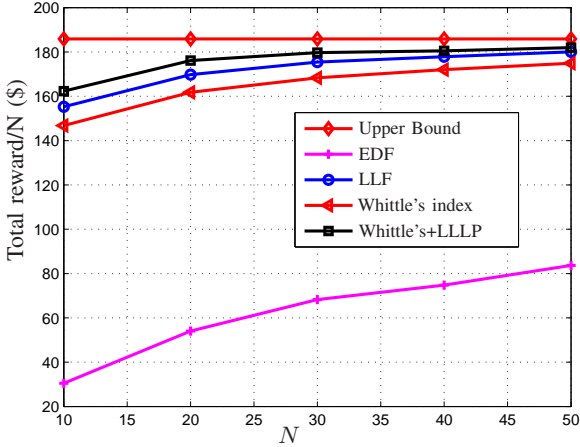}
  \caption{Performance comparison with constant processing cost: ${c[t]=0.5}$, ${Q(0,0)=0.3}$, ${\bar{T}=12}$, ${\bar{B}=9}$, ${\beta=0.999}$, ${F(B)=0.2B^2}$, ${M/N=0.5}$.}
  \label{fig:constantCost}
\end{figure}

In Figure \ref{fig:constantCost}, we considered the case when ${M/N=0.5}$ and varied the maximum queue length $N$. We observed that the Whittle's index policy with LLLP interchange and LLF achieved similar performance, since both policies roughly followed the least laxity first principle. The performance of these two policies was close to the performance upper bound.
The EDF policy performed poorly because it did not take the remaining job length into account.
The gap between the Whittle's index policy and the Whittle's index policy with LLLP enhancement came from the reordering of jobs with positive laxity (cf. the discussion following
Definition \ref{def:WhittleIndexPolicy}).

\subsection{Dynamic Processing Cost}

For the dynamic processing cost case, we used the real-time electricity price signal from the California Independent System Operator (CAISO) and trained a Markovian model that described the
marginal processing costs (cf. Sections III and V of \cite{Kwon&Xu&Gautam:TSG15}). Each time slot of the
constructed Markov chain (on processing cost) lasted for 1 hour.
For each time slot, the real-time price was quantized into discrete price states, and the transition probability (of the Markov chain) was simply the frequency the price changes from one state to another.

\begin{figure}[ht]
\centering
\psfrag{Num of arm}{{\small $M/N$}}
\psfrag{Reward per arm}{{\small Total reward/N (\$)}}
\psfrag{UB}{{\scriptsize Upper Bound}}
\psfrag{EDF}{{\scriptsize EDF}}
\psfrag{LLF}{{\scriptsize LLF}}
\psfrag{Valley filling}{{\scriptsize Valley filling}}
\psfrag{Whittls index}{{\scriptsize Whittle's index}}
\psfrag{LLLLLLLLLLLLP-Index}{{\scriptsize Whittle's+LLLP}}
  \includegraphics[width=.48\textwidth]{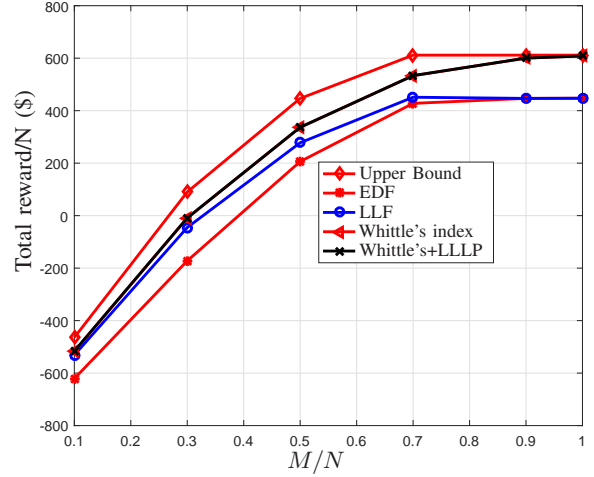}
  \caption{Performance comparison with dynamic processing cost: ${Q(0,0)=0.3}$, ${\bar{T}=12}$, ${\bar{B}=9}$, ${\beta=0.999}$, ${F(B)=0.2B^2}$, ${N=10}$.}
  \label{fig:varyRatio}
\end{figure}

In Figure~\ref{fig:varyRatio}, we fixed the job arrival process and the maximum queue length ${N=10}$ and varied the processing capacity $M$. When the processing limit was low and $M/N$ was small, there were not enough processors to finish all jobs,
 and the non-completion penalty dominated the processing profit.
 In this case, the performance of different policies was close due to the little flexibility constrained by the limited processing resource.
 When the processing capacity was adequate and ${M/N=1}$, all jobs could be finished on time. In this case, the Whittle's index policy solved the problem optimally and achieved the upper bound (which was in correspondence with Proposition \ref{prop:optimal}).
 The LLLP interchange never happened because the Whittle's index policy followed
 the LLLP principle in this case. EDF and LLF did not utilize any information about the stochastic processing cost process and achieved sub-optimal performance.
 When the processing capacity constraint was neither too tight ($M/N\approx0$) nor too loose ($M/N\approx1$), the LLLP principle tended to break large unfinished jobs (with long remaining processing time) into smaller jobs and therefore improved the overall performance by
 processing more tasks when processing cost was low and reducing the non-completion penalty.

\begin{figure}
\centering
\psfrag{Num of arm}{{\small $N$}}
\psfrag{Reward per arm}{{\small Total reward (\$)}}
\psfrag{UB}{{\scriptsize Upper Bound}}
\psfrag{EDF}{{\scriptsize EDF}}
\psfrag{LLF}{{\scriptsize LLF}}
\psfrag{SL}{{\scriptsize Valley filling}}
\psfrag{Index-Idle}{{\scriptsize Whittle's index}}
\psfrag{Index-Idle-NoBack}{{\scriptsize Whittle's+LLLP}}
\psfrag{Whittls index}{{\scriptsize Whittle's index}}
\psfrag{LLLLLLLLLLLLP-Index}{{\scriptsize Whittle's+LLLP}}
\psfrag{Valley filling}{{\scriptsize Valley filling}}
  \includegraphics[width=.48\textwidth]{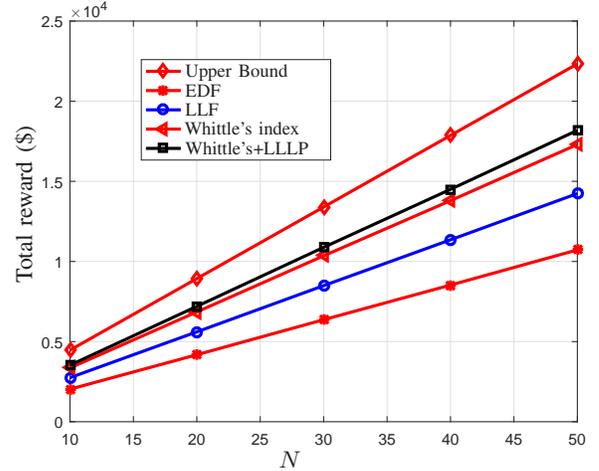}
  \caption{Performance comparison with dynamic processing cost: ${Q(0,0)=0.3}$, ${\bar{T}=12}$, ${\bar{B}=9}$, ${\beta=0.999}$, ${F(B)=0.2B^2}$, ${M/N=0.5}$.}
  \label{fig:dynamicCost}
\end{figure}

In Figure~\ref{fig:dynamicCost}, we compared the performance of different policies by fixing ratio $M/N=0.5$ and varying the maximum queue length $N$.
Both the EDF and LLF policies sought to activate as many jobs as possible, up to the processing capacity  $M$.
The Whittle's index policy, on the other hand, took pricing fluctuation into account: it processed more jobs at price valley and keept processors idle when the processing cost was high. Based on the Whittle's index policy, the LLLP enhancement further reduced the penalty of unfinished jobs and improved the performance of the Whittle's index policy.
The total reward achieved by the Whittle's index with LLLP enhancement policy was more than 1.7 times of that obtained by EDF, and the performance gap between
the Whittle's index with LLLP policy and the LLF policy was over $25\%$. We also noticed that the LLLP principle improved the Whittle's index policy by around $10\%$.

\subsection{Asymptotic Optimality}
\begin{figure}[ht]
\centering
\psfrag{M/N}{{\scriptsize $M$}}
\psfrag{Whittls index}{{\scriptsize $G(s)-G_{\mbox{W}}(s)$}}
\psfrag{EDF}{{\scriptsize $G(s)-G_{\tiny \mbox{EDF}}(s)$}}
\psfrag{LLF}{{\scriptsize $G(s)-G_{\tiny \mbox{LLF}}(s)$}}
\psfrag{Reward per arm}{{\scriptsize Gap-to-optimality}}
\psfrag{LLLLLLLLLLLLLLLLP-Index}{{\scriptsize $G(s)-G_{\mbox{W+LLLP}}(s)$}}
\psfrag{Upper Bound}{{\scriptsize Bound on gap in (\ref{eqn:Poisson})}}
  \includegraphics[width=.52\textwidth]{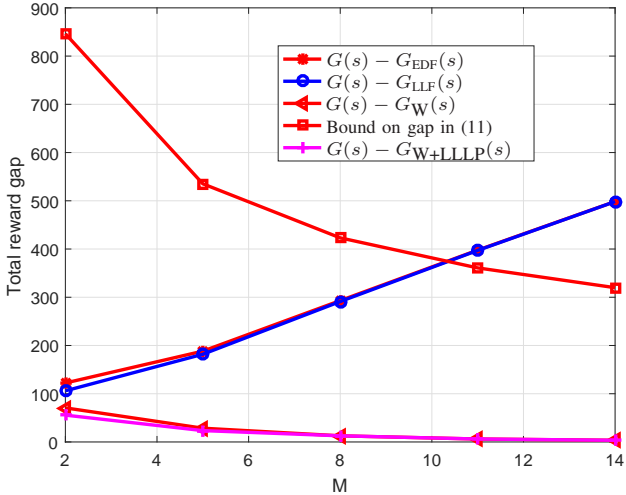}
  \caption{Gap-to-optimality of three different index policies under dynamic processing cost: ${Q(0,0)=0.3}$, ${\bar{T}=12}$, ${\bar{B}=9}$, ${\beta=0.999}$, ${F(B)=0.2B^2}$, $\mu=M$, $N=1000$. }
  \label{fig:asy}
\end{figure}
In Figure~\ref{fig:asy},  simulation results are presented to compare the performance achieved by various heuristic policies
and to validate the theoretic results established in Lemma \ref{lem:gap}.

 In this simulation, we fixed the queue size ${N=1000}$ and varied the processing capacity $M$ as a parameter. The arrival sequence within $\bar{T}$ time slots was generated from a Poisson process with mean ${\mu=M}$.
The dynamic cost evolved according to a Markovian model
that was trained using real-time
electricity price signals from CAISO. Each time slot of the constructed Markov chain
 lasted for 1 hour, and the entire simulation horizon lasted for 300 days (with $24\times300$ time slots).

The EDF and LLF policies did not take into account the dynamics of processing costs,
and their gap-to-optimality increased as both the job arrival rate and processing capacity grew as shown in Figure \ref{fig:asy}.
On the other hand, the gap between the total rewards achieved by the Whittle's index policy and the optimal policy quickly decreased to zero as $M$ increases.
We note that the Whittle's index policy's actual gap-to-optimality
was less than the performance gap bound derived in  (\ref{eqn:Poisson}), as shown in
Theorem \ref{thm:Poisson}. We also showed in Figure \ref{fig:asy} the gap-to-optimality of the LLLP enhanced Whittle's index policy.
 The performance gap of the Whittle's index policy and the LLLP enhanced one was small because the arrival traffic was relatively light.

\subsection{Hard Deadlines}\label{sec:hardDeadlines}
In this subsection, we examine the performance of the proposed algorithms in a setting with hard deadlines. In this setting, we seek to
finish as many jobs as possible regardless of the processing cost. Our framework can incorporate the hard deadline scenario by
 setting the non-completion penalty much higher than processing costs.
In our simulation,
 we set the processing cost ${c=0.95}$ and considered a linear penalty function with a slope of $10$, $F(B)=10B$.
 In this setting, it was optimal ({\it i.e.,} reward maximizing) to finish as many jobs as possible. 

\begin{figure}[ht]
\centering
\psfrag{Num of arm}{{\scriptsize $M/N$}}
\psfrag{Finished Job}{{\scriptsize Completion ratio}}
\psfrag{Whittls index}{{\scriptsize Whittle's index}}
\psfrag{LLLLLLLLLLLLP-Index}{{\scriptsize Whittle's+LLLP}}
\psfrag{IGLLSP}{{\scriptsize Whittle's+LLSP}}
\psfrag{EDF}{{\scriptsize EDF}}
\psfrag{LLF}{{\scriptsize LLF}}
  \includegraphics[width=.5\textwidth]{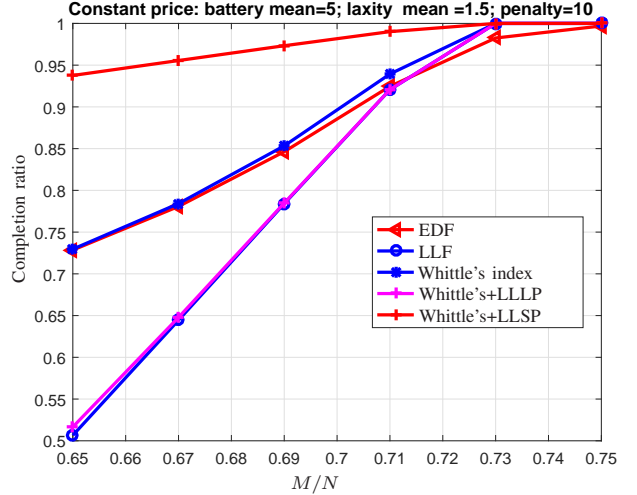}
  \caption{Job completion ratio: ${c[t]=0.95}$, ${Q(0,0)=0.3}$, ${\bar{T}=12}$, ${\bar{B}=9}$, ${\beta=0.999}$, ${F(B)=10B}$, $N=100$.}
  \label{fig:hardDL}
\end{figure}

The ratios of completed jobs achieved by various algorithms are plotted in Figure \ref{fig:hardDL}.
We noted that the Whittle's index policy outperformed the EDF and LLF policies.
Although the LLLP principle improved the Whittle's index policy in the sense of total reward, it completed fewer jobs as LLLP can result in many small unfinished jobs.
Interestingly, we observed from Figure \ref{fig:hardDL} that the Less Laxity Smaller Processing time (LLSP) principle could significantly enhance the job completion ratio
achieved by the Whittle's index policy.
The LLSP interchange is the same as the LLLP interchange (introduced in Section \ref{sec:op4}), except that priority will be given to smaller unfinished jobs instead of larger unfinished ones.

\subsection{Validation of Assumption A5-A6}\label{sec:validationOfA6}
We conducted numerical experiments to evaluate the effect of the independent arrival assumption in A5-A6. We considered two scenarios.
In the first scenario, the job arriving at each position followed an independent binomial distribution (according to A5-A6).
In the second scenario, the aggregate job arriving at the system followed a Poisson distribution with the same mean as that in the first scenario. When a job arriveed at the system, it was randomly assigned to one of the empty positions following a uniform distribution.
\begin{figure}[ht]
\centering
\psfrag{N}{{\scriptsize Slots in the queue $N$}}
\psfrag{Reward per arm}{{\scriptsize Total Reward}}
\psfrag{Whittls index}{{\scriptsize Whittle's index w. A5-A6.}}
\psfrag{UB}{{\scriptsize Upper bound w. A5-A6.}}
\psfrag{Whittles+LLLP}{{\scriptsize Whittle's+LLLP w. A5-A6.}}
\psfrag{EDF}{{\scriptsize EDF w. A5-A6.}}
\psfrag{LLF}{{\scriptsize LLF w. A5-A6.}}
\psfrag{Whittles index Poisson LLLLLLLLLL}{{\scriptsize Whittle's index w. Poisson arr.}}
\psfrag{Upper bound w. Poisson Arrivals.}{{\scriptsize Upper bound w. Poisson arr.}}
\psfrag{Whittles+LLLLLLLLP Poisson}{{\scriptsize Whittle's+LLLP w. Poisson arr.}}
\psfrag{EDF Poisson}{{\scriptsize EDF w. Poisson arr.}}
\psfrag{LLF Poisson}{{\scriptsize LLF w. Poisson arr.}}
  \includegraphics[width=.5\textwidth]{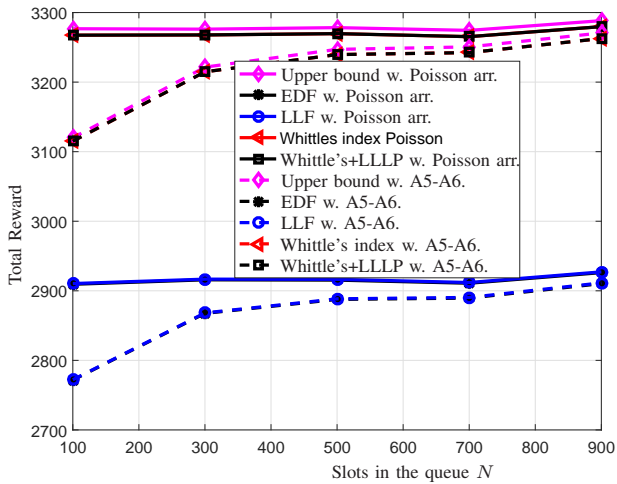}
  \caption{Comparison between Poisson arrival and independent arrival under dynamic processing cost: ${Q(0,0)=0.3}$, ${\bar{T}=12}$, ${\bar{B}=9}$, ${\beta=0.999}$, ${F(B)=0.2B^2}$, $M=10$, $\mu=M$. }
  \label{fig:PoissonArrival}
\end{figure}

We let the number of available processors ${M=10}$ and fixed the mean of the total job arrivals (within $\bar{T}$ time slots) as ${\mu=M}$. Shown in Figure \ref{fig:PoissonArrival}, as the number of available positions in the queue increased, the performance of different algorithms under A6 converged to its counterpart under Poisson arrival. 

%% file: conclusion_v1.tex
We consider the problem of large scale deadline scheduling---a problem that has a wide range of applications in calling centers, cloud computing, and EV charging.  In such settings, it is essential to develop highly efficient and online scheduling algorithms.  To this end, the index policy considered in this paper are attractive for its implementation simplicity and versatility in incorporating various operation uncertainties. It is particularly reassuring that the upper bound on the gap-to-optimality of the Whittle's index policy converges to zero, thus establishing the asymptotic optimality of the Whittle's index policy in the light traffic regime.

%% file: proof_v11.tex
\section{Proof of Theorem \ref{thm:indexability}}\label{proof:indexability}
In \cite{Graczova&Jacko:14OR}, the indexability of the bi-dimension state model is proved without arrivals. In this appendix, we provide an elementary proof for the indexability of the RMAB problem
formulated in (\ref{eqn:MAB}) with random arrivals. In particular, we will show that for any state $\tilde{s}$ of an arm,
there is a critical $\nu(\tilde{s})$ such that if and only if $\nu\ge\nu(\tilde{s})$ the first term in the Bellman equation (\ref{eqn:BellmanNuSubsidy}) is larger than or equal to the second term in a single arm $\nu$-subsidy problem.

\subsection{Indexability of Dummy Arms}\label{sec:secAppDummy}
The indexability of dummy arms is straightforward. For ${i\in\{N+1,\cdots,N+M\}}$, there is no job arrival, and only the processing cost evolves.
The Bellman equation of the $\nu$-subsidy problem is given by
\[\begin{array}{l}
V_i^\nu(0,0,c_j)=\max\{\beta\sum_kP_{j,k}V^\nu_i(0,0,c_k)+\nu,\\[4pt]
\quad \mathrel{\phantom{V_i^\nu(0,0,c_j)=\max\{}}\beta\sum_kP_{j,k}V^\nu_i(0,0,c_k)\}.
\end{array}\]

If and only if ${\nu\ge0}$, the first term is larger than the second term and it is optimal to deactivate the dummy arm. Otherwise, the active action is optimal. So a dummy arm is indexable and its Whittle's index is ${\nu_i(0,0,c_j)=0}$.

\subsection{Indexability of Regular Arms}

\begin{proof}
We now prove the indexability of regular arms by induction.
We first show that the Whittle's index $\nu_i(T,B,c_j)$ exists for ${T\le1}$ and all $B$ and $c_j$,
and establish some useful properties
for
the difference of the value function
${g^\nu(T,B,c_j)\triangleq V_i^\nu(T,B+1,c_j)-V_i^\nu(T,B,c_j)}$
     for the case with ${T=1}$. Then, under the conditions that the Whittle's index $\nu_i(T,B,c_j)$ exists and the property of $g^\nu(T,B,c_j)$ holds for ${T=t-1}$, we show $\nu_i(T,B,c_j)$ exists, and the property of $g^\nu(T,B,c_j)$ holds for ${T=t}$.

\subsubsection{${T=0}$}
 There is no job waiting in the position. The Bellman equation is stated as
    \[V_i^\nu(0,0,c_j)=\max\{\nu+\beta W_{j}^\nu,\beta W_{j}^\nu\},\]
        where
            \[
            \begin{array}{ll} \\
               W_{j}^\nu&=\sum_{T'}\sum_{B}\sum_{k}Q(T',B)P_{j,k}V_i^\nu(T',B,c_k)\\
                &\mathrel{\phantom{=}}+Q(0,0)\sum_{k}P_{j,k}V_i^\nu(0,0,c_k)
                        \end{array}
                            \] is the expected reward of possible arrivals. If and only if ${\nu\ge0}$, the first term is larger and the passive action is optimal. Thus ${\nu_i(0,0,c_j)=0}$.
\subsubsection{$T=1$}
 There are two cases.
    \bitem
    \item If $B=0$, the Bellman equation is stated as
\[V_i^\nu(1,0,c_j)=\max\{\nu+\beta W_{j}^\nu,\beta W_{j}^\nu\}.\]
Thus ${\nu_i(1,0,c_j)=0}$.
    \item If $B\ge1$, the Bellman equation is stated as
    \[\hspace{-1em}\begin{array}{l}
    V_i^\nu(1,B,c_j)=\max\{\nu-F(B)+\beta W_{j}^\nu,\\
    \mathrel{\phantom{V_i^\nu(1,B,c_j)=\max\{}}1-c_j-F(B-1)+\beta W_{j}^\nu\}.
    \end{array}\]
    \normalsize
   If and only if ${\nu\ge1-c_j+F(B)-F(B-1)}$, the passive action is optimal.
\eitem

   Thus the Whittle's index for ${T=1}$ exists, and the closed-form is given by
    \beq\label{eqn:closed-formT1}
    \begin{array}{l}
    \nu_i(1,B,c_j)\\
    ~~~~=\left\{
    \begin{array}{ll}
    0,&\mbox{if~} B=0;\\
    1-c_j+F(B)-F(B-1),&\mbox{if~} B\ge1.
    \end{array}
    \right.\end{array}
    \eeq

    Let the difference of the value function be
     \[
     g^\nu(T,B,c_j)\triangleq V_i^\nu(T,B+1,c_j)-V_i^\nu(T,B,c_j).
       \]
       We note that the difference of the value function is continuous and piecewise linear in $\nu$. Specially, denote $\mathcal{G}$ as a set of functions of $\nu$ such that ${g(\nu)\in\mathcal{G}}$ if and only if $g(\nu)$ is a continuous piecewise linear function in $\nu$, there exist $\underline{\nu}$ and $\bar{\nu}$ such that, ${\partial g(\nu)/\partial\nu\ge-1}$ when ${\nu\in[\underline{\nu},\bar{\nu}]}$, and ${\partial g(\nu)/\partial\nu=0}$ when ${\nu\notin[\underline{\nu},\bar{\nu}]}$. We show that, when $T=1$, $g^\nu(T,B,c_j)\in\mathcal{G}$.

       \bitem
\item If $B=0$, \[g^\nu(1,B,c_j)=V_i^\nu(1,1,c_j)-V_i^\nu(1,0,c_j).\]
\bitem
\item If $\nu_i(1,1,c_j)>\nu_i(1,0,c_j)=0$,
\[\hspace{-1em}
g^\nu(1,B,c_j)
=\left\{
\begin{array}{ll}
1-c_j,&\mbox{if~}\nu<0;\\
1-c_j-\nu,&\mbox{if~}0\le\nu<\nu_i(1,1,c_j);\\
-F(1),&\mbox{if~}\nu_i(1,1,c_j)\le\nu.
\end{array}
\right.\]

\item If ${\nu_i(1,1,c_j)\le\nu_i(1,0,c_j)=0}$,
\[\hspace{-1em}
g^\nu(1,B,c_j)
=\left\{
\begin{array}{ll}
1-c_j,&\mbox{if~}\nu<\nu_i(1,1,c_j);\\
\nu-F(1),&\mbox{if~}\nu_i(1,1,c_j)\le\nu<0;\\
-F(1),&\mbox{if~}0\le\nu.
\end{array}
\right.\]
\eitem

\item If $B\ge1$, \[{g^\nu(1,B+1,c_j)=V_i^\nu(1,B+1,c_j)-V_i^\nu(1,B,c_j)}.\]
Since ${\nu_i(1,B+1,c_j)\ge\nu_i(1,B,c_j)}$ by (\ref{eqn:closed-formT1}),
\[\hspace{-1em}g^\nu(1,B,c_j)
=\left\{
\begin{array}{l}
F(B-1)-F(B),\\
\mbox{if~}\nu<\nu_i(1,B,c_j);\\
1-c_j-\nu,\\
\mbox{if~}\nu_i(1,B,c_j)\le\nu<\nu_i(1,B+1,c_j);\\
F(B)-F(B+1),\\
\mbox{if~}\nu_i(1,B+1,c_j)\le\nu.
\end{array}
\right.\]
So $g^\nu(1,B,c_j)$ is continuous piecewise linear in $\nu$, and there exist ${\underline{\nu}}$ and $\bar{\nu}$ such that ${\partial g^\nu(1,B,c_j)/\partial \nu\ge-1}$ when ${\nu\in[\underline{\nu},\bar{\nu}]}$ and ${\partial g^\nu(1,B,c_j)/\partial \nu=0}$ otherwise.
\eitem

\subsubsection{$T\ge2$} Assuming the Whittle's index $\nu_i(T,B,c_j)$ exits and $g^\nu(T,B,c_j)\in\mathcal{G}$ for ${T=t-1}$, we show $\nu_i(T,B,c_j)$ exits and $g^\nu(T,B,c_j)\in\mathcal{G}$ for the case $T=t$.

First, existence of ${\nu_i(T,B,c_j)}$ when $T=t$.
\bitem
\item If $B=0$, the Bellman equation is stated as follows.
\[\hspace{-1em}\begin{array}{l}
    V_i^\nu(t,0,c_j)
    =\max\{\beta\sum_kP_{j,k}V_i^\nu(t-1,0,c_k)+\nu,\\
    \mathrel{\phantom{V_i^\nu(t,0,c_j)
    =\max\{}}\beta\sum_kP_{j,k}V_i^\nu(t-1,0,c_k)\}.
    \end{array}\]
    If and only if ${\nu\ge0}$, the first term is larger than the second term and the passive action is optimal. Thus ${\nu_i(t,0,c_j)=0}$.

\item If $B\ge1$, the Bellman equation is stated as follows.
    \begin{equation}\label{eqn:BLargerThanOen}
    \hspace{-1em}
    \begin{array}{l}
    V_i^\nu(t,B,c_j)\\
    ~~~~=\max\{\beta\sum_kP_{j,k}V_i^\nu(t-1,B,c_k)+\nu,\\
    ~~~~\mathrel{\phantom{=\max\{}}\beta\sum_kP_{j,k}V_i^\nu(t-1,B-1,c_k)+1-c_j\}.
    \end{array}
    \end{equation}
    Denote the difference between the two actions as
    \[
    \begin{array}{ll}
     f^\nu(t,B,c_j)&\triangleq\beta\sum_kP_{j,k}g^\nu(t-1,B-1,c_k)\\
   &\mathrel{\phantom{\triangleq}}+\nu-(1-c_j),
    \end{array}
    \]
    where \[
    \begin{array}{l}
    g^\nu(t-1,B-1,c_k)\\
    ~~~~=V_i^\nu(t-1,B,c_k)-V_i^\nu(t-1,B-1,c_k).
    \end{array}\]

Since ${g^\nu(t-1,B-1,c_k)\in\mathcal{G}}$ by assumption, ${f^\nu(t,B,c_j)}$ is continuous and piece-wise linear in $\nu$. Let
\[\begin{array}{l}
\underline{\nu}(t,B,c_j)\triangleq\min_k\underline{\nu}(t-1,B-1,c_k),\\
\bar{\nu}(t,B,c_j)\triangleq\max_k\bar{\nu}(t-1,B-1,c_k),
\end{array}\]
where ${\partial g^\nu(t-1,B-1,c_k)/\partial \nu\ge-1}$ if and only if ${\nu\in[\underline{\nu}(t-1,B-1,c_k),\bar{\nu}(t-1,B-1,c_k)]}$.
We have
\[\begin{array}{l}
\partial f^\nu(t,B,c_j)/\partial\nu\\
~~~~=\left\{
\begin{array}{ll}
\ge0,&\mbox{if~}\nu\in[\underline{\nu}(t,B,c_j),\bar{\nu}(t,B,c_j)]; \\
1,& \mbox{otherwise.}
\end{array}
\right.
\end{array}\]

So $f^\nu(t,B,c_j)$ is continuous and non-decreasing in $\nu$. When ${\nu=-\infty}$, ${f^\nu(t,B,c_j)=-\infty}$. When ${\nu=+\infty}$, ${f^\nu(t,B,c_j)=+\infty}$. Thus there is a cross point of $f^\nu(t,B,c_j)$ and the $\nu$-axis. Define ${\nu_i(t,B,c_j)\triangleq\min_{\nu}\{f^\nu(t,B,c_j)=0\}}$. If and only if ${\nu\ge\nu_i(t,B,c_j)}$, the first term in (\ref{eqn:BLargerThanOen}) is larger or equal to the second term and the passive action is optimal. By definition, $\nu_i(t,B,c_j)$ is the Whittle's index.
\eitem
The existence of $\nu_i(t,B,c_j)$ is shown.

Next we show $g^\nu(t,B,c_j)\in\mathcal{G}$.

\bitem
\item If $B=0$, \[g^\nu(t,B,c_j)=V_i^\nu(t,1,c_j)-V_i^\nu(t,0,c_j).\]
\bitem
\item If $\nu_i(t,1,c_j)>\nu_i(t,0,c_j)=0$,
\[\hspace{-3em}
\begin{array}{l}
g^\nu(t,0,c_j)\\
~~~~=\left\{
\begin{array}{ll}
1-c_j,&
\mbox{if~}\nu<0;\\
1-c_j-\nu,&
 \mbox{if~}0\le\nu<\nu_i(t,1,c_j);\\
\beta\sum_kP_{j,k}g^\nu(t-1,0,c_k),&
 \mbox{if~}\nu_i(t,1,c_j)\le\nu.
\end{array}
\right.
\end{array}
\]
\item If $\nu_i(t,1,c_j)\le\nu_i(t,0,c_j)=0$,
\[\hspace{-1em}
g^\nu(t,0,c_j)=\left\{
\begin{array}{l}
1-c_j,\\
\mbox{if~}\nu<\nu_i(t,1,c_j);\\
\nu+\beta\sum_kP_{j,k}g^\nu(t-1,0,c_k),\\
\mbox{if~}\nu_i(t,1,c_j)\le\nu<0;\\
\beta\sum_kP_{j,k}g^\nu(t-1,0,c_k),\\
\mbox{if~}0\le\nu.
\end{array}
\right.
\]
\eitem
Thus, ${g^\nu(t,0,c_j)}$ is a linear combination of ${g^\nu(t-1,0,c_k)}$. Since ${g^\nu(t-1,0,c_k)\in\mathcal{G}}$ for all $c_k$ by assumption, we have ${g^\nu(t,0,c_j)}\in\mathcal{G}$ as well.
\item If $B\ge1$, \[g^\nu(t,B,c_j)=V_i^\nu(t,B+1,c_j)-V_i^\nu(t,B,c_j).\]
\bitem
\item If $\nu_i(t,B+1,c_j)>\nu_i(t,B,c_j)$,
\[\hspace{-1em}
\begin{array}{l}
g^\nu(t,B,c_j)\\
~~~~=\left\{
\begin{array}{l}
\beta\sum_kP_{j,k}g^\nu(t-1,B-1,c_k),\\
\mbox{if~}\nu<\nu_i(t,B,c_j);\\
1-c_j-\nu,\\
\mbox{if~}\nu_i(t,B,c_j)\le\nu<\nu_i(t,B+1,c_j);\\
\beta\sum_kP_{j,k}g^\nu(t-1,B,c_k),\\
\mbox{if~}\nu_i(t,B+1,c_j)\le\nu.
\end{array}
\right.
\end{array}
\]
\item If $\nu_i(t,B+1,c_j)\le\nu_i(t,B,c_j)$,
\[\hspace{-3em}
\begin{array}{l}
g^\nu(t,B,c_j)\\
~~~~=\left\{
\begin{array}{l}
\beta\sum_kP_{j,k}g^\nu(t-1,B-1,c_k),\\
\mbox{if~}\nu<\nu_i(t,B+1,c_j);\\
\beta\sum_kP_{j,k}[g^\nu(t-1,B,c_k)+g^\nu(t-1,B-1,c_k)]\\
+\nu-(1-c_j),\\
\mbox{if~}\nu_i(t,B+1,c_j)\le\nu<\nu_i(t,B,c_j);\\
\beta\sum_kP_{j,k}g^\nu(t-1,B,c_k),\\
\mbox{if~}\nu_i(t,B,c_j)\le\nu.
\end{array}
\right.
\end{array}
\]
\eitem
\eitem

Clearly, ${g^\nu(t,B,c_j)}$ is a linear combination of ${g^\nu(t-1,B,c_k)}$ and ${g^\nu(t-1,B-1,c_k)}$. Since ${g^\nu(t-1,B,c_k)\in\mathcal{G}}$ for all $B$ and $c_k$ by assumption, we have ${g^\nu(t,B,c_j)}\in\mathcal{G}$ as well.

Thus, by induction, the Whittle's index $\nu_i(T,B,c_j)$ exists and ${g^\nu(T,B,c_j)\in\mathcal{G}}$ for all $T,B$, and $c_j$.

\end{proof}

\section{Proof of Theorem \ref{thm:closed-form}}\label{proof:closed-form}
\begin{proof}
Since the processing cost $c_0$ is constant, we omit the cost in the state of arms for simplicity. In Appendix \ref{sec:secAppDummy} we have shown that the Whittle's index of the dummy arms is $\nu_i(0,0)=0$.
For regular arms, we showed in (\ref{eqn:closed-formT1}) that ${\nu_i(1,0)=0}$ and ${\nu_i(1,B)=1-c_0+F(B)-F(B-1)}$ when $B\ge1$. Next, we show the closed-form of the Whittle's index for the case of ${T\ge2}$ using induction.

\subsection {$T=2$}
The discussion is divided into two conditions.
\bitem
\item If $B=1$,
 \[
 \begin{array}{l}
V_i^\nu(2,1)=\max\{\nu+\beta V_i^\nu(1,1),\\
 \mathrel{\phantom{V_i^\nu(2,1)=\max\{}}1-c_0+\beta V_i^\nu(1,0)\}.
 \end{array}
 \]
 The difference between active and passive actions
\[\hspace{-1em}
\begin{array}{l}
f^\nu(2,1)\\
~~~~=\nu-(1-c_0)+\beta g^\nu(1,0)\\
~~~~=\left\{
\begin{array}{ll}
\nu-(1-\beta)(1-c_0),& \mbox{if~}\nu<0;\\
(1-\beta)[\nu-(1-c_0)],& \mbox{if~}0\le\nu<1-c_0+F(1);\\
\nu-(1-c_0)-\beta F(1),& \mbox{if~}1-c_0+F(1)\le\nu;\\
\end{array}
\right.
\end{array}
\]
\normalsize
equals $0$ when ${\nu=1-c_0}$. Thus ${\nu_i(2,1)=1-c_0}$.

\item If $B\ge2$, the Bellman equation is stated as follows.
\[\hspace{-1em}
\begin{array}{l}
V_i^\nu(2,B)=\max\{\nu+\beta V_i^\nu(1,B),\\
\mathrel{\phantom{V_i^\nu(2,B)=\max\{}}1-c_0+\beta V_i^\nu(1,B-1)\}.
\end{array}
\]
Let ${\Delta F(B)=F(B)-F(B-1)}$. The difference between active and passive actions
\[\hspace{-1em}
\begin{array}{l}
f^\nu(2,B)\\
~~~~=\nu-(1-c_0)+\beta g^\nu(1,B-1)\\
~~~~=\left\{
\begin{array}{l}
\nu-(1-c_0)-\beta\Delta F(B-1),\\
\mbox{if~}\nu<1-c_0+\Delta F(B-1);\\
(1-\beta)[\nu-(1-c_0)],\\
\mbox{if~}1-c_0+\Delta F(B-1)\le\nu<1-c_0+\Delta F(B);\\
\nu-(1-c_0)+\beta\Delta F(B),\\
\mbox{if~}1-c_0+\Delta F(B)\le\nu;\\
\end{array}
\right.
\end{array}
\]
equals $0$ when ${\nu=1-c_0+\beta[F(B-1)-F(B-2)]}$. Thus ${\nu_i(2,B)=1-c_0+\beta[F(B-1)-F(B-2)]}$ when $B\ge2$.

So (\ref{eqn:closedForm}) is true when $T=2$.
\eitem

\subsection {${T>2}$}
Assume equation (\ref{eqn:closedForm}) holds when ${T=t-1}$, we show that it holds when ${T=t}$.
\bitem
\item If $B=1$,
\[\begin{array}{l}
V_i^\nu(t,B)=\max\{\nu+\beta V_i^\nu(t-1,1),\\
\mathrel{\phantom{V_i^\nu(t,B)=\max\{}}1-c_0+\beta V_i^\nu(t-1,0)\}.
\end{array}
\]
The difference between actions is
\[\hspace{-1em}
\begin{array}{l}
f^\nu(t,1)\\
~~~~=\nu-(1-c_0)+\beta g^\nu(t-1,0)\\
~~~~=\left\{
\begin{array}{ll}
\nu-(1-\beta)(1-c_0),& \mbox{if~}\nu<0;\\
(1-\beta)[\nu-(1-c_0)],& \mbox{if~}0\le\nu<1-c_0;\\
\nu-(1-c_0)+\beta^2g^\nu(t-2,0),& \mbox{if~}1-c_0\le\nu.\\
\end{array}
\right.
\end{array}
\]
The last case can be rewritten as
\[
\begin{array}{l}
\nu-(1-c_0)+\beta^2g^\nu(t-2,0)\\
~~~~=(1-\beta)[\nu-(1-c_0)]+\beta[\nu-(1-c_0)]\\
~~~~\mathrel{\phantom{=}}+\beta^2[V_i^\nu(t-2,1)-V_i^\nu(t-2,0)],
\end{array}\]
which equals $0$ when ${\nu=1-c_0}$ since by assumption ${\nu_i(t-1,1)=1-c_0}$.  Thus ${\nu_i(t,1)=1-c_0}$.

\item If $2\le B\le t-2$, the difference between actions is stated as follows.
\[
\begin{array}{l}
f^\nu(t,B)\\
~~~~=\nu-(1-c_0)+\beta g^\nu(t-1,B-1)\\
~~~~=\left\{
\begin{array}{l}
\beta^2g^\nu(t-2,B-2)
+\nu-(1-c_0),\\[2pt] \mbox{if~}\nu<1-c_0;\\[3pt]
\beta^2g^\nu(t-2,B-1)
+\nu-(1-c_0),\\[2pt]  \mbox{if~}1-c_0\le\nu.
\end{array}
\right.
\end{array}
\]
The latter case equals $0$ when ${\nu=1-c_0}$ because ${\nu_i(t-1,B)=1-c_0}$ when ${2\le B\le t-2}$ by assumption. Thus ${\nu_i(t,B)=1-c_0}$ when ${2\le B\le t-2}$.

\item If $B=t-1$,
\[\hspace{-1em}\begin{array}{ll}
f^\nu(t,B)&=\nu-(1-c_0)+\beta g^\nu(t-1,B-1)\\
&=\left\{
\begin{array}{l}
\nu-(1-c_0)+\beta^2g^\nu(t-2,B-2),\\[2pt]
 \mbox{if~}\nu<1-c_0;\\[3pt]
(1-\beta)[\nu-(1-c_0)],\\[2pt]
  \mbox{if~}1-c_0\le\nu<1-c_0+\beta^{t-2}F(1);\\[3pt]
\nu-(1-c_0)+\beta^2g^\nu(t-2,B-1),\\[2pt]
 \mbox{if~}1-c_0+\beta^{t-2}F(1)\le\nu;
\end{array}
\right.
\end{array}
\]
equals $0$ when $\nu=1-c_0$. So ${\nu_i(t,B)=1-c_0}$ when ${B=t-1}$.

\item If ${B\ge t}$,
\begin{equation}\label{eqn:differenceOfValueFunction}
\hspace{-1em}
\begin{array}{ll}
f^\nu(t,B)&=\nu-(1-c_0)+\beta g^\nu(t-1,B-1)\\
&=\left\{
\begin{array}{l}
\nu-(1-c_0)+\beta^2g^\nu(t-2,B-2),\\[2pt]
 \mbox{if~}\nu<\nu_i(t-1,B-1);\\[2pt]
(1-\beta)[\nu-(1-c_0)],\\[2pt]
 \mbox{if~}\nu_i(t-1,B-1)\le\nu<\nu_i(t-1,B);\\[2pt]
\nu-(1-c_0)+\beta^2g^\nu(t-2,B-1),\\[2pt]
 \mbox{if~}\nu_i(t-1,B)\le\nu.\\
\end{array}
\right.
\end{array}
\end{equation}
If $\nu<\nu_i(t-1,B-1)$, according to (\ref{eqn:closedForm})
\[
\begin{array}{ll}
\nu&<\nu_i(t-1-T',B-1-T')\\
&\le\nu_i(t-1-T',B-T'),
\end{array}
\]
for all ${0\le T'\le t-1}$. Thus the first case of (\ref{eqn:differenceOfValueFunction}) can be written as
\[\hspace{-1em}
\begin{array}{l}
\nu-(1-c_0)+\beta^2g^\nu(t-2,B-2)\\
~~~~=\nu-(1-c_0)+\beta^3 g^\nu(t-3,B-3)\\
~~~~=\cdots\\
~~~~=\nu-(1-c_0)+\beta^{t-1} g^\nu(1,B-t+1)\\
~~~~=\nu-(1-c_0)+\beta^{t-1}[-F(B-t+1)+F(B-t)].
\end{array}
\]
\eitem
As a result, when ${\nu=1-c_0+\beta^{t-1}[F(B-t+1)-F(B-t)]}$, the first case in equation (\ref{eqn:differenceOfValueFunction}) equals $0$ . Thus when $B\ge t$, the closed-form of index is stated as:
\[
\nu_i(t,B)=1-c_0+\beta^{t-1}[F(B-t+1)-F(B-t)].
\]
We therefore conclude that (\ref{eqn:closedForm}) holds when ${T=t}$. By induction, we have established (\ref{eqn:closedForm}) for all $T$.

\end{proof}

\section{Proof of Proposition \ref{prop:optimal}}\label{proof:optimalityOfWhittleIndexPolicy}
\begin{proof}
In this appendix, we show that, when ${M=N}$, the Whittle's index policy optimally solves the RMAB problem defined in (\ref{eqn:MAB}), which is equivalent to the MDP problem in (\ref{eqn:originPro}). We first prove that the Whittle's index policy optimally solves the reward maximizing problem of a single arm without constraint in Lemma \ref{col:optimalityOfWhittleindex}. Then we show that the Whittle's index policy is optimal for the RMAB problem defined in (\ref{eqn:MAB}) when ${M=N=1}$. In the end, we show the optimality of the Whittle's index policy for (\ref{eqn:MAB})  when ${M=N}$ in general.

\subsection{Optimality in single arm problem}
First, we claim that the Whittle's index policy optimally solves the single arm problem with dynamic cost and no constraint. The extended state ${\tilde{s}=(T,B,c)}$ includes the state of the job and the processing cost. The Bellman equation of the single arm problem is given by:
\begin{equation}\label{eqn:singleArmMAB}
V_i(\tilde{s})=\max\{R_0(\tilde{s})+\beta(\mathcal{L}_{0}V_i)(\tilde{s}), R_1(\tilde{s})+\beta(\mathcal{L}_{1}V_i)(\tilde{s})\},
\end{equation}
 where action ${a=1}$ means to activate the arm and $a=0$ means to leave it passive.

The Whittle's index is defined by introducing a  {\it $\nu$-subsidy problem},
which is a modified version of the single arm problem defined in (\ref{eqn:singleArmMAB}).
In the $\nu$-subsidy problem, whenever the passive action is taken, the scheduler receives an extra reward $\nu$  \cite{Whittle:1988JAP}. The single arm problem defined in (\ref{eqn:singleArmMAB}) is simply the case when the subsidy $\nu=0$.

The Bellman equation for the $\nu$-subsidy problem is given by
\beq\label{eqn:nu-problem}
V_i^\nu(\tilde{s})=\max\{R_0(\tilde{s})+\nu+\beta(\mathcal{L}_{0}V_i^\nu)(\tilde{s}), R_1(\tilde{s})+\beta(\mathcal{L}_{1}V_i^\nu)(\tilde{s})\},
\eeq
where $V_i^\nu$ is the value function for the $\nu$-subsidy problem.

Now define a Whittle's index policy $\pi_1$ for a single arm (either regular or dummy arm) $\nu$-subsidy problem as to activates the arm if and only if ${\nu_i(\tilde{s})>\nu}$. Thus we have the optimality of $\pi_1$ as follows.

\begin{lem}\label{col:optimalityOfWhittleindex}
The Whittle's index policy $\pi_1$ is optimal for the single arm $\nu$-subsidy problem defined in (\ref{eqn:nu-problem}). In particular, when ${\nu=0}$, $\pi_1$ is optimal for the single arm problem defined in (\ref{eqn:singleArmMAB}).
\end{lem}

\begin{proof}
We have shown in Appendix \ref{proof:indexability} that
the Whittle's index defined in Definition \ref{def:W} exists,
and therefore the Whittle's index policy $\pi_1$ is well defined.
By Definition \ref{def:W}, for any state $\tilde{s}$ such that ${\nu_i(\tilde{s})>\nu}$, the first term in (\ref{eqn:nu-problem}) is strictly smaller than the second term. The Whittle's index policy $\pi_1$ activates the arm  and obtains the second term as an expected reward which satisfies the Bellman equation in this case.

For ${\nu=\nu_i(\tilde{s})}$, the first term is greater or equal to the second term in the Bellman equation by Definition \ref{def:W}. The indexability result proved in Appendix \ref{proof:indexability} guarantees that the passive set grows monotonously which implies that this inequality is true for any ${\nu\ge\nu(\tilde{s})}$. Thus, for any state $\tilde{s}$ such that ${\nu(\tilde{s})\le\nu}$, the Whittle's index policy $\pi_1$ leaves the arm passive and obtains the first term as the expected reward, satisfying the Bellman equation.

Thus, $\pi_1$ satisfies the Bellman equation (\ref{eqn:nu-problem}) and is therefore optimal for the single arm $\nu$-subsidy problem. In particular, when ${\nu=0}$, $\pi_1$ is optimal for the single arm problem and satisfies the Bellman equation (\ref{eqn:singleArmMAB}).
\end{proof}

\subsection{Optimality when ${M=N=1}$}
Now we consider the problem (\ref{eqn:MAB}) with ${M=N=1}$:
 we have a regular arm and a dummy arm, and at each time, we are required to activate exact one arm.
For this constrained two-arm problem, the state is defined as ${(\tilde{s},\textbf{0})=(T,B,c,0,0)}$, where $\tilde{s}$ is the extended state of the regular arm and $\textbf{0}=(0,0)$ the state of the dummy arm. The action ${a'=1}$ means to activate the regular arm, and ${a'=0}$ represents activating the dummy arm.

The state of the dummy arm will always be $\textbf{0}$. The dummy arm yields no reward regardless of the taken action. Thus the state transition of two-arm problem
is equivalent to the state transition in problem (\ref{eqn:singleArmMAB}), \ie ${P((\tilde{s},\textbf{0}),(\tilde{s}',\textbf{0})|a')=P(\tilde{s},\tilde{s}'|a)}$.
The rewards of the two-arm problem can be presented by the rewards of
 the single arm problem  in (\ref{eqn:singleArmMAB}):
\[\begin{array}{l}
R'_1(\tilde{s},\textbf{0})=R_1(\tilde{s}),\\
R'_0(\tilde{s},\textbf{0})=R_0(\tilde{s}).
\end{array}
\]

The Whittle's index policy for the two-arm problem (denoted by $\pi_2$) activates the regular arm when ${(\nu(\tilde{s})>\nu(\textbf{0})=0)}$, and activates the dummy arm (leaving the regular arm passive) otherwise.

 When $\pi_1$ faces state $\tilde{s}$ and $\pi_2$ faces state $(\tilde{s},\textbf{0})$ for the same realization $\tilde{s}$, the actions of two policies are the same.  $\pi_2$ will activate the regular arm in the two-arm problem if and only if $\pi_1$ activates the arm in the single arm problem, and vice versa. Since the reward, transition and the action of these two policies are the same, the value functions will be the same. Denoting the value function of $\pi_1$ and $\pi_2$ by $V_{\pi_1}(\tilde{s})$ and $H_{\pi_2}(\tilde{s},\textbf{0})$, we have $H_{\pi_2}(\tilde{s},\textbf{0})=V_{\pi_1}(\tilde{s})$. Since  $V_{\pi_1}(\tilde{s})$ satisfies the Bellman equation (\ref{eqn:singleArmMAB}), we have
\[\begin{array}{l}
\mathrel{\phantom{=}}H_{\pi_2}(\tilde{s},\textbf{0})\\
=\max\{R_0(\tilde{s})+\beta(\mathcal{L}_{0}H_{\pi_2})(\tilde{s},\textbf{0}), R_1(\tilde{s})+\beta(\mathcal{L}_{1}H_{\pi_2})(\tilde{s},\textbf{0})\}\\
=\max\{R'_0(\tilde{s},\textbf{0})+\beta(\mathcal{L}_{0}H_{\pi_2})(\tilde{s},\textbf{0}),\\ \mathrel{\phantom{=\max\{}}R'_1(\tilde{s},\textbf{0})+\beta(\mathcal{L}_{1}H_{\pi_2})(\tilde{s},\textbf{0})\},
\end{array}\]
which is in fact the Bellman equation for the constrained two-arm problem.
The Whittle's index policy satisfies the Bellman equation for the two-arm problem and is therefore optimal.

\subsection{Optimality when ${M=N}$}
Finally, we argue that the Whittle's index policy is optimal for the multi-arm problem defined in (\ref{eqn:MAB}) when ${M=N}$.
We have $N$ regular arms and $N$ dummy arms. At each time, we activate exact $N$ arms. We can pair each regular arm with a dummy arm and implement the Whittle's index policy for each pair. The action of each regular arm is decoupled, and the total reward is simply the sum of reward from all the $N$ regular arms. The Whittle's index policy optimally  solves the problem of each pair, and is therefore optimal for the original problem in (\ref{eqn:MAB}).

\end{proof}

\section{Proof of Lemma~\ref{lem:gap}}\label{proof:asymptotic}
\begin{proof}
For Problem (\ref{eqn:MAB}), 
we use $\pi_{\tiny\mbox{W}}$ to denote the Whittle's index policy with the processing limit $M$ (that activates the $M$ arms with highest indices at each time). Specially, when the limit is loose, \eg ${M=N}$ denote the Whittle's index policy (that activates the $N$ arms with highest indices at each time) by $\pi_{\tiny \mbox{N}}$.

 In Appendix \ref{proof:optimalityOfWhittleIndexPolicy}, we have shown that when ${M=N}$, the Whittle's index policy is optimal. Thus the reward of $\pi_{\tiny \mbox{N}}$ serves as an upper bound of the optimal reward for any case with ${M\le N}$, \ie
\[G^N_{\pi_{\tiny \mbox{N}}}(s)\ge G^N(s)\ge G^N_{\tiny\mbox{W}}(s),\]
where $G^N_{\pi_{\tiny \mbox{N}}}(s)$ is the reward collected from $\pi_{\tiny \mbox{N}}$, $G^N_{\tiny\mbox{W}}(s)$ the reward collected from the Whittle's index policy $\pi_{\tiny\mbox{W}}$ when $M\le N$, and $G^N(s)$ the maximum reward defined in (\ref{eqn:originPro}).

In this appendix, we establish an upper bound of the difference of the value functions of $\pi_{\tiny\mbox{W}}$ and $\pi_{\tiny \mbox{N}}$, ${G^N_{\pi_{\tiny \mbox{N}}}(s)-G^N_{\tiny\mbox{W}}(s)}$, which serves as an upper bound of the gap-to-optimality of the Whittle's index policy, ${G^N(s)-G^N_{{\tiny\mbox{W}}}(s)}$. We first quantify ${G^N_{\pi_{\tiny \mbox{N}}}(s)-G^N_{{\tiny\mbox{W}}}(s)}$ by the number of different actions in the processing sequences resulted by $\pi_{\tiny\mbox{W}}$ and $\pi_{\tiny \mbox{N}}$. Then we relate the number of different actions to the number of job arrivals in Lemma \ref{lem:boundA}, which gives us the result in (\ref{eqn:bound}).

Note that due to the lack of capacity limit, policy $\pi_{\tiny \mbox{N}}$ activates a regular arm if and only if its Whittle's index is positive. On the other hand, the policy $\pi_{\tiny\mbox{W}}$  activates a regular arm if and only if its index belongs to the largest $M$ positive ones. Due to the capacity limit $M$, if facing the same trajectory of processing cost and the same sequence of arrivals, the two policies $\pi_{\tiny \mbox{N}}$ and $\pi_{\tiny\mbox{W}}$ will generate different processing sequences on a single job. As shown in Figure \ref{fig:actionDiff}, the processing sequences of a job $J_i$ with arrival time $r$ and departure time $d$ determined by  $\pi_{\tiny \mbox{N}}$ and $\pi_{\tiny\mbox{W}}$ are plotted. We define two events as follows.
 \bitem
 \item Event $\textbf{A}$: $\pi_{\tiny \mbox{N}}$ processes $J_i$ but $\pi_{\tiny\mbox{W}}$ does not.
\item Event $\textbf{B}$: $\pi_{\tiny\mbox{W}}$ processes $J_i$ but $\pi_{\tiny \mbox{N}}$ does not.
 \eitem

\begin{figure}[h]
\centering
\psfrag{N}{{\small $\pi_{\tiny \mbox{N}}$}}
\psfrag{A}{{\small $\textbf{A}$}}
\psfrag{B}{{\small $\textbf{B}$}}
\psfrag{M}{{\small $\pi_{\tiny\mbox{W}}$}}
\psfrag{r}{{\small $r$}}
\psfrag{d}{{\small $d$}}
  \includegraphics[width=.45\textwidth]{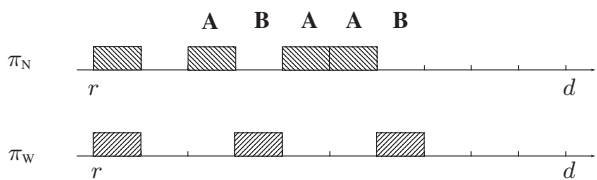}
  \caption{Processing sequences generated by $\pi_{\tiny \mbox{N}}$ and $\pi_{\tiny\mbox{W}}$ on a single job.}
  \label{fig:actionDiff}
\end{figure}

\bitem
\item If event $\textbf{A}$ happens at time $t$, the instant reward difference between $\pi_{\tiny \mbox{N}}$ and $\pi_{\tiny\mbox{W}}$ is bounded, \ie
\[
R_i^{\pi_{\tiny \mbox{N}}}[t]-R_i^{\pi_{\tiny\mbox{W}}}[t]\le |1-c_{\min}|,
\]
where $R_i^{\pi}[t]$ is the instant reward collected from $J_i$ by policy $\pi$ at time $t$.

\item If event $\textbf{B}$ happens at time $t$, the instant reward difference between $\pi_{\tiny \mbox{N}}$ and $\pi_{\tiny\mbox{W}}$ is also bounded, \ie
\[
R_i^{\pi_{\tiny \mbox{N}}}[t]-R_i^{\pi_{\tiny\mbox{W}}}[t]\le |1-c_{\max}|.
\]

\item At the deadline of $J_i$, the difference of unfinished job length resulting from two policies is bounded by the number of event $\textbf{A}$. Thus the penalty difference of two policies is also bounded, \ie

\[
\begin{array}{ll}
    F_i^{\pi_{\tiny \mbox{N}}}[d]-F_i^{\pi_{\tiny\mbox{W}}}[d]&\le F(B+\sum_{t=r}^d\mathds{1}({\textbf{A}[t]}))-F(B)\\
    &\le F(\bar{B})\sum_{t=r}^d\mathds{1}({\textbf{A}[t]}),
\end{array}
\]
where $F_i^{\pi}[d]$ is the penalty of $J_i$ resulted by $\pi$ at deadline $d$, $B$ is the left over job size under $\pi_{\tiny \mbox{N}}$ of job $i$, $\mathds{1}({\textbf{A}[t]})=1$ if and only if event $\textbf{A}$ happens at time $t$, and $F(\bar{B})$ is the maximum penalty that can incur to a job.
\eitem
The reward difference collected from $J_i$ up to time $t<d$ is the sum of the first two cases, \ie
\[
 \begin{array}{l}
\sum_{h=r}^t\beta^h(R_i^{\pi_{\tiny \mbox{N}}}[h]-R_i^{\pi_{\tiny\mbox{W}}}[h])\\
~~~~\le|1-c_{\min}|\sum_{h=r}^t\mathds{1}({\textbf{A}[h]})\beta^h\\
~~~~~~~+|1-c_{\max}|\sum_{h=r}^t\mathds{1}({\textbf{B}[h]})\beta^h.
 \end{array}
 \]
\normalsize
The difference up to deadline $t=d$ is the sum of the three cases, \ie
\[
 \begin{array}{l}
\sum_{h=r}^d\beta^h(R_i^{\pi_{\tiny \mbox{N}}}[h]-R_i^{\pi_{\tiny\mbox{W}}}[h])\\
~~~~\le|1-c_{\min}|\sum_{h=r}^d\mathds{1}({\textbf{A}[h]})\beta^h\\
\mathrel{\phantom{~~~~\le}}+|1-c_{\max}|\sum_{h=r}^d\mathds{1}({\textbf{B}[h]})\beta^h\\
\mathrel{\phantom{~~~~\le}}+F(\bar{B})\beta^d\sum_{h=r}^d\mathds{1}({\textbf{A}[h]}).
 \end{array}
 \]
For each time $t\in[r,d]$, we enlarge the penalty term and get a general bound as follows.
\begin{equation}\label{eqn:penaltyGeneralBound}
 \begin{array}{l}
\sum_{h=r}^t\beta^h(R_i^{\pi_{\tiny \mbox{N}}}[h]-R_i^{\pi_{\tiny\mbox{W}}}[h])\\
~~~~\le |1-c_{\min}|\sum_{h=r}^t\mathds{1}({\textbf{A}[h]})\beta^h\\
\mathrel{\phantom{~~~~\le }}+|1-c_{\max}|\sum_{h=r}^t\mathds{1}({\textbf{B}[h]})\beta^h\\
\mathrel{\phantom{~~~~\le }}+F(\bar{B})\sum_{h=r}^t\mathds{1}({\textbf{A}[h]})\beta^h.
 \end{array}
 \end{equation}

Note that the cumulative number of event $\textbf{A}$ happened up to any fixed time $t$ is always larger than the number of event $\textbf{B}$. Formally, we state the following lemma to illustrate the relationship between event $\textbf{A}$ and $\textbf{B}$. The proof is delayed to Appendix \ref{proof:boundB}.
\begin{lem}\label{lem:boundB}
Denote $\mathds{1}({\textbf{A}[t]})$ as whether event $\textbf{A}$ happens at $t$.
Denote $\#\textbf{A}[t]$ as the cumulative number of event $\textbf{A}$ happened from $r$ to time $t$. Define $\mathds{1}({\textbf{B}[t]})$ and $\#\textbf{B}[t]$ respectively. For any $t\in[r,d]$,
\[
\begin{array}{c}
\#\textbf{A}[t]=\sum_{h=r}^t\mathds{1}({\textbf{A}[h]})\ge\#\textbf{B}[t]=\sum_{h=r}^t\mathds{1}({\textbf{B}[h]}),\\
\sum_{h=r}^t\mathds{1}({\textbf{A}[h]})\beta^h\ge\sum_{h=r}^t\mathds{1}({\textbf{B}[h]})\beta^h.
\end{array}\]
\end{lem}
So the reward difference in (\ref{eqn:penaltyGeneralBound}) is bounded as follows.
\[
 \begin{array}{l}
\sum_{h=r}^t\beta^h(R_i^{\pi_{\tiny \mbox{N}}}[h]-R_i^{\pi_{\tiny\mbox{W}}}[h])\\
~~~~~\le (|1-c_{\min}|+F(\bar{B})+|1-c_{\max}|)\sum_{h=r}^t\mathds{1}({\textbf{A}[h]})\beta^h.
 \end{array}
 \]

Now we want to quantify the cumulative number of event $\textbf{A}$. Event $\textbf{A}$ happens only when there are more than $M$ jobs with positive Whittle's index in the system under $\pi_{\tiny\mbox{W}}$. This event can only occur when there are at least $M$ jobs in the queue. To bound the number of event $\textbf{A}$, we have the following lemma. The proof is delayed to Appendix \ref{proof:boundA}.

\begin{lem}\label{lem:boundA}
Let $I^N[t]$ be the number of jobs admitted to the system within ${[t-\bar{T}+1,t]}$. Then for any $t$, \[{\mathds{1}}({\textbf{A}[t]})\le \mathds{1}({I^N[t]>M}).\]
\end{lem}

Thus for each job, we have
\begin{equation}\label{eqn:singleEV}
 \sum_{h=r}^t\beta^h(R_i^{\pi_{\tiny \mbox{N}}}[h]-R_i^{\pi_{\tiny\mbox{W}}}[h])\le C\sum_{h=r}^t\beta^h\mathds{1}({I^N[h]>M}),
 \end{equation}
for any $t$, where $C=(|1-c_{\min}|+F(\bar{B})+|1-c_{\max}|)$.

If we sum arrivals and take expectation, we have the difference of expected value function bounded as follows.
\[\begin{array}{ll}
G^N_{\pi_{\tiny \mbox{N}}}(s)-G^N_{\tiny\mbox{W}}(s)&\le C\sum_t\beta^t\mathbb{E}\big[\mathds{1}({I^N[t]>M})I^N[t]\big]\\[3pt]
&=C\mathbb{E}\big[\mathds{1}({I^N[t]>M})I^N[t]\big]/(1-\beta).
\end{array}
\]

Since $G^N_{\pi_{\tiny \mbox{N}}}(s)$ is an upper bound of $G^N(s)$, we have
\[\begin{array}{ll}
G^N(s)-G^N_{\tiny\mbox{W}}(s)&\le \frac{C}{1-\beta}\mathbb{E}\big[I^N[t]\mathds{1}(I^N[t]>M)\big]\\
&=\frac{C}{1-\beta}\mathbb{E}\big[I^N[t]|I^N[t]>M\big]\Pr(I^N[t]>M),
\end{array}
\]
which is the expression (\ref{eqn:bound}) in Lemna \ref{lem:gap}.

\end{proof}

\subsection{Proof of Lemma \ref{lem:boundB}}\label{proof:boundB}
\begin{proof}
At time $t$, we denote the remaining job size of $J_i$ under policy $\pi_{\tiny \mbox{N}}$ and $\pi_{\tiny\mbox{W}}$ by $B_{\tiny \mbox{N}}[t]$ and
$B_{\tiny\mbox{W}}[t]$, respectively.
When ${B_{\tiny\mbox{W}}[t]=B_{\tiny \mbox{N}}[t]}$, job $J_i$ has the same state and Whittle's index under both polices. If $\pi_{\tiny\mbox{W}}$ processes $J_i$, which means the Whittle's index of $J_i$ is positive, $\pi_{\tiny \mbox{N}}$ also processes $J_i$. Since at the arrival $r$, $B_{\tiny\mbox{W}}[r]=B_{\tiny \mbox{N}}[r]$, event $\textbf{B}$ can only happen when $B_{\tiny\mbox{W}}[t]>B_{\tiny \mbox{N}}[t]$, which means event $\textbf{A}$ must have happened before.

This also implies that $B_{\tiny\mbox{W}}[t]\ge B_{\tiny \mbox{N}}[t]$ for all $t$.
\end{proof}

\subsection{Proof of Lemma~\ref{lem:boundA}} \label{proof:boundA}
\begin{proof}
Recall that the remaining job size under $\pi_{\tiny\mbox{W}}$ is always larger than the one under $\pi_{\tiny \mbox{N}}$, \ie  $B_{\tiny\mbox{W}}[t]\ge B_{\tiny \mbox{N}}[t]$. Whenever $\pi_{\tiny \mbox{N}}$ processes some job $J_i$, the Whittle's index of this job under $\pi_{\tiny \mbox{N}}$ must be positive. If we can show that the Whittle's index is monotonically increasing in $B$, the index under policy $\pi_{\tiny\mbox{W}}$ must also be positive, and $\pi_{\tiny\mbox{W}}$ will also process this job if the capacity limit allows. Thus that event $\textbf{A}$ happens must imply that there are more than $M$ jobs with positive Whittle's index, which requires the number of admitted jobs larger than $M$, \ie $I^N[t]>M$.

In this subsection, we show that the Whittle's index is indeed increasing in $B$ when the index is positive and the value function is concave when ${\nu>0}$ by induction. That is, ${\nu_i(T,B+1,c_j)\ge\nu_i(T,B,c_j)}$,  if ${\nu_i(T,B,c_j)>0}$ and ${V_i^{\nu}(T,B,c_j)}$ is concave when ${\nu>0}$.
\subsubsection {${T=1}$} The Whittle's index is
 \[\hspace{-1em}   \nu_i(1,B,c_j)=\left\{
    \begin{array}{ll}
    0,&\mbox{if~} B=0;\\
    1-c_j+F(B)-F(B-1),&\mbox{if~} B\ge1.
    \end{array}
    \right.
    \]
    If $\nu_i(1,B,c_j)>0$, $\nu_i(1,B+1,c_j)>\nu_i(1,B,c_j)$ due to the convexity of $F(B)$.

    The value function is concave in $B$ when $\nu>0$.
   \[\hspace{-1em}
   \begin{array}{l}
   V_i^\nu(1,B+2,c_j)-2V_i^\nu(1,B+1,c_j)+V_i^\nu(1,B,c_j)\\
   ~~~~=\left\{
   \begin{array}{l}
   -F(B+2)+2F(B+1)-F(B),\\
    \mbox{if } \nu_i(1,B,c_j)<\nu, \nu_i(1,B+1,c_j)<\nu,\\
    \mbox{and } \nu_i(1,B+2,c_j)<\nu;\\
    1-c_j-\nu+F(B+1)-F(B),\\
    \mbox{if } \nu_i(1,B,c_j)<\nu, \nu_i(1,B+1,c_j)<\nu,\\
    \mbox{and }\nu\le  \nu_i(1,B+2,c_j);\\
    \nu-1+c_j+F(B)-F(B+1),\\
    \mbox{if } \nu_i(1,B,c_j)<\nu\le\nu_i(1,B+1,c_j);\\
    -F(B+1)+2F(B)-F(B-1),\\
    \mbox{if } \nu\le\nu_i(1,B,c_j);\\
   \end{array}
   \right.
   \\
   ~~~~\le0
   \end{array}
   \]
   The first and last cases are negative because of convexity of the penalty. The second and third cases are negative because of the definition of $\nu_i(1,B,c_j)$.
\subsubsection{${T>1}$} Assume ${\nu_i(T,B+1,c_j)\ge\nu_i(T,B,c_j)}$ when $\nu_i(T,B,c_j)>0$, and $V_i^\nu(T,B,c_j)$ is concave in $B$ when $\nu>0$ for ${T=t-1}$. We show that these properties are true for ${T=t}$.

The difference of the activate and deactivate actions at state ${(t,B+1,c_j)}$ is given by
   \[
   \begin{array}{l}
   f^\nu(t,B+1,c_j)\\
   ~~~~=\beta\sum P_{j,k}[V_i^{\nu}(t-1,B+1,c_k)-V_i^{\nu}(t-1,B,c_k)]\\
   ~~~~\mathrel{\phantom{=}}+\nu-1+c_j\\
   ~~~~=   \beta\sum P_{j,k}[V_i^{\nu}(t-1,B+1,c_k)-2V_i^{\nu}(t-1,B,c_k)\\
    ~~~~\mathrel{\phantom{=}}+V_i^{\nu}(t-1,B-1,c_k)]\\

  ~~~~\mathrel{\phantom{=}}+\beta\sum P_{j,k}[V_i^{\nu}(t-1,B,c_k)-V_i^{\nu}(t-1,B-1,c_k)]\\
   ~~~~\mathrel{\phantom{=}}+\nu-1+c_j\\
     ~~~~=\beta\sum P_{j,k}[V_i^{\nu}(t-1,B+1,c_k)-2V_i^{\nu}(t-1,B,c_k)\\
    ~~~~\mathrel{\phantom{=}}+V_i^{\nu}(t-1,B-1,c_k)]\\
   ~~~~\mathrel{\phantom{=}}+f^\nu(t,B,c_j).
   \end{array}
   \]
   When ${\nu=\nu_i(t,B,c_j)>0}$, we have ${f^\nu(t,B,c_j)=0}$ according to the definition of ${\nu_i(t,B,c_j)}$. The first term in the above equation is negative due to the concavity of the value function when ${\nu>0}$. We thus have ${f^\nu(t,B+1,c_j)\le0}$ when ${\nu=\nu_i(t,B,c_j)>0}$, which implies ${\nu_i(t,B+1,c_j)\ge\nu_i(t,B,c_j)}$.

We have shown the monotonicity of the Whittle's index when ${T=t}$. Next we show the concavity of the value functions for $T=t$ when $\nu>0$.
    \[\hspace{-1em}
   \begin{array}{l}
   V_i^\nu(t,B+2,c_j)+V_i^\nu(t,B,c_j)-2V_i^\nu(t,B+1,c_j)\\
   ~~~~=\left\{
   \begin{array}{l}
   \beta\sum P_{j,k}V_i^\nu(t-1,B+2,c_k)\\
   -2\beta\sum P_{j,k}V_i^\nu(t-1,B+1,c_k)\\
   +\beta\sum P_{j,k}V_i^\nu(t-1,B,c_k),\\
    \mbox{if } \nu_i(t,B,c_j)<\nu, \nu_i(t,B+1,c_j)<\nu,\\
     \mbox{and } \nu_i(t,B+2,c_j)<\nu;\\
    \beta\sum P_{j,k}V_i^\nu(t-1,B,c_k)\\
    -\beta\sum P_{j,k}V_i^\nu(t-1,B+1,c_k)\\
    +1-c_j-\nu,\\
    \mbox{if } \nu_i(t,B,c_j)<\nu, \nu_i(t,B+1,c_j)<\nu,\\
    \mbox{and }\nu\le  \nu_i(t,B+2,c_j);\\
    \beta\sum P_{j,k}V_i^\nu(t-1,B+1,c_k)\\
    -\beta\sum P_{j,k}V_i^\nu(t-1,B,c_k)\\
    +\nu-(1-c_j),\\
    \mbox{if } \nu_i(t,B,c_j)<\nu\le\nu_i(t,B+1,c_j);\\
    \beta\sum P_{j,k}V_i^\nu(t-1,B+1,c_k)\\
    -\beta\sum P_{j,k}2V_i^\nu(t-1,B,c_k)\\
    +\beta\sum P_{j,k}V_i^\nu(t-1,B-1,c_k),\\
    \mbox{if } \nu\le\nu_i(t,B,c_j);
   \end{array}
   \right.
   \\
   ~~~~\le0.
   \end{array}
   \]
The first and fourth terms are less than zero because by the assumption the value function is concave when ${\nu>0}$ for ${t-1}$. The second and third terms are negative because of the definition of $\nu_i(t,B+1,c_j)$. So the value function $V^\nu_i(t,B,c_j)$ is concave in $B$ when $\nu>0$.

By induction, we have ${\nu_i(T,B+1,c_j)\ge\nu_i(T,B,c_j)}$ when ${\nu_i(T,B,c_j)>0}$, and $V_i^\nu(T,B,c_j)$ is concave in $B$ when $\nu>0$ for all $T$.
\end{proof}

\section{Proof of Corollary \ref{col:infiniteN}}\label{proof:infiniteN}
\begin{proof}
Consider a job $J$ arrive at time $t$. When it arrives at the queue, there are at most $I[t]$ jobs waiting in the system and the probability of job $J$ gets rejected is no more than $I[t]/N$. Since $\mathbb{E}(I[t])<\infty$, for any $\varepsilon>0$, there exists  $K$ such that $Pr(I[t]>K)<\varepsilon/2$. So for ${N>\max\{K(K+1)/\varepsilon,K\}}$, the rejection probability of $J$ is less than $\varepsilon$,
\[
\begin{array}{l}
\mathrel{\phantom{\le}}Pr(J \mbox{ is rejected})\\
\le \sum_{i=1}^N Pr(I[t]=i)i/N+\sum_{N+1}^\infty Pr(I[t]=i) \\
\le \sum_{i=1}^K Pr(I[t]=i)i/N+Pr(I[t]>K)\\
\le K(K+1)/2N+Pr(I[t]>K)\\
\le \varepsilon
\end{array}
\]
The probability that no job get rejected is no smaller than $(1-\epsilon)^{I[t]}$ and $I^N[t]$ converges to $I[t]$ in probability as follows.
\[
Pr(I^N[t]=I[t])\ge (1-\epsilon)^{I[t]}\rightarrow1, \mbox{ as } N\rightarrow\infty.
\]

\end{proof}

\section{Proof of Theorem \ref{thm:Poisson}}\label{proof:Poisson}
\begin{proof}
For a Poisson process $I[t]$ with mean $\mu$, we have the expression as follows.
\[
\mathbb{E}[I[t]|I[t]>M]Pr(I[t]>M)=\mu Pr(I[t]\ge M)
\]
For any $M>\mu-1$, we have the inequality as follows \cite{Klar:2000PEIS}.
\beq\label{eqn:poissonDis}
\begin{array}{ll}
\mu Pr(I[t]\ge M)&<\mu Pr(I[t]=M)/(1-\dfrac{\mu}{M+1})\\[5pt]
&=\mu^{M+1}e^{-\mu}(M+1)/{[(M+1-\mu)M!]}\\[5pt]
&\le\dfrac{\mu^{M+1}e^{M-\mu}(M+1)}{\sqrt{2\pi}M^{M+1/2}(M+1-\mu)}\\[5pt]
&=O(\dfrac{\mu e^{-\mu}}{\sqrt{M}})
\end{array}
\eeq
where the second inequality is because of Stirling formula. When ${\mu\le M/e}$, the right-hand side decreases to zero, which indicates the asymptotic optimality of the Whittle's index.
\end{proof}